\documentclass[11pt,letterpaper]{article}

\pdfoutput=1

       \usepackage{amsfonts}
       \usepackage{stmaryrd}
       \usepackage{latexsym,amssymb,mathrsfs,fancyhdr}
       \usepackage{amsmath,amssymb}
       \usepackage{tikz}
       \usetikzlibrary{matrix,arrows}
\usepackage{extarrows}
\usepackage{todonotes}

       \font\tenmsb=msbm10
       \font\sevenmsb=msbm7
       \font\fivemsb=msbm5
       \catcode`\@=11
       \ifx\amstexloaded@\relax\catcode`\@=\active
       \endinput\else\let\amstexloaded@\relax\fi
       \def\spaces@{\space\space\space\space\space}
       \def\spaces@@{\spaces@\spaces@\spaces@\spaces@\spaces@}
       \def\space@.  {\futurelet\space@\relax}
       \space@.   %
       \def\Err@#1{\errhelp\defaulthelp@\errmessage{AmS-TeX error: #1}}
       \def\relaxnext@{\let\next\relax}
       \def\accentfam@{7}
       \def\noaccents@{\def\accentfam@{0}}
       \def\Cal{\relaxnext@\ifmmode\let\next\Cal@\else
       \def\next{\Err@{Use \string\Cal\space only in math mode}}\fi\next}
       \def\Cal@#1{{\Cal@@{#1}}}
       \def\Cal@@#1{\noaccents@\fam\tw@#1}
       \def\Bbb{\relaxnext@\ifmmode\let\next\Bbb@\else
       \def\next{\Err@{Use \string\Bbb\space only in math mode}}\fi\next}
       \def\Bbb@#1{{\Bbb@@{#1}}}
       \def\Bbb@@#1{\noaccents@\fam\msbfam#1}
       \newfam\msbfam
       \textfont\msbfam=\tenmsb
       \scriptfont\msbfam=\sevenmsb
       \scriptscriptfont\msbfam=\fivemsb

\def\N{\mathbb{N}}
\def\Z{\mathbb{Z}}
\def\R{\mathbb{R}}

\def\C{\mathbb{C}}

\def\im{ \mathrm{Im}\,}
\def\supp{ \mathrm{supp}\,}
\def\e{\epsilon}
 \def\ep{\varepsilon}
\def\be{\varepsilon}

 \def\pd{ \frac{\partial}{\partial \nu}}
 \def\lag{  \langle}
 \def\rag{ \rangle}
\def\dt{\nabla \wedge}
\def\nt{\nu \wedge}

\def\exi{e^{-\frac{i}{h}x\cdot \xi}}
\newcommand{\su}{\subseteq}

\newcommand{\ra}{{\rightarrow}}
\newcommand{\ma}{{\mapsto}}
 \def\om{\omega}
\def\Om{\Omega}
\def\pd {\partial}
\def\pn {\partial_{\nu}}
\def\Ga{\Gamma}

\def\skipaline{\removelastskip\vskip12pt plus 1pt minus 1pt}

\def\Proof{\removelastskip\skipaline \noindent \it Proof.  \rm}

\newtheorem{ Remark}{ Remark}[section]

\newcommand{\qed}{\nobreak \hfill \ifvmode \relax \else
      \vrule height0.75em width0.5em depth0.25em \fi}

\newcommand{\beq}{\begin{equation} }
\newcommand{\eeq}{\end{equation} }

\usepackage[german,english]{babel}

\newtheorem{theo}{Theorem}[section]

\newtheorem{rem}{Remark}[section]
\newtheorem{prop}{Proposition}[section]
\newtheorem{lem}{Lemma}[section]
\newtheorem{iteration lemma}{iteration Lemma}[section]

\begin{document}

\setlength{\columnsep}{5pt}
\title{Inverse problems for nonlinear Helmholtz Schr\"odinger equations and time-harmonic Maxwell's equations with partial data}
\author{~Xuezhu  Lu}

\maketitle

\tableofcontents

\begin{abstract}
We consider Calder\'{o}n's  inverse boundary value problems   for a class of nonlinear Helmholtz Schr\"{o}dinger equations and Maxwell's equations in a bounded domain in $\R^n$. The main method is the higher order linearization of the Dirichlet-to-Neumann map of the corresponding equations. The local uniqueness of the   linearized partial data Calder\'{o}n's  inverse problem is obtained following  \cite{DKSU}. The Runge approximation properties and unique continuation principle allow to extend to global situations. Simultaneous recovery of some unknown cavity$/$boundary and coefficients are given as some applications.

\end{abstract}

\section{Introduction and Main Results}

  In this paper, we  study the partial data inverse boundary value problems for a class of nonlinear Helmholtz Schr\"{o}dinger equations and Maxwell's equations in a {bounded} domain.

 We first consider the inverse boundary value problems of the following nonlinear Helmholtz Schr\"odinger equation
\begin{align}\label{main1}
\left\{\begin{array}{ll}
  -\Delta u(x)-k^2 u(x)+q(x)u(x)^2 =0~~&\text{in~}\Omega,\\
   u(x) =f(x)~~&\text{on~} \pd \Om,
   \end{array}
   \right.
\end{align}
where $f\in C^s(\partial\Omega)$, $s>1$ with $s\not\in\Z$. Here,   $\Omega$ is a bounded open set in $\R^n$ with $n\geq 2.$ Assume that the boundary  $\partial \Omega$ is smooth.

In Appendix \ref{sec:appdx} we show that except for a discrete set of real-valued $k$, the Dirichlet problem (\ref{main1}) has a unique small solution $u\in W^{2,p}(\Omega)$ for sufficiently small boundary data $f\in W^{2-\frac1p, p}(\pd \Om)$, for $p>n$. It is also shown in \cite{LLLS1} that   the boundary value problem admits a unique solution in $C^s(\overline{\Omega})$ when $\|f\|_{C^s(\partial\Omega)}$ is small in the case of $k=0$; and then  the full boundary Dirichlet-to-Neumann (DN) map  $\Lambda_q$ is well defined as follows:
 $$\Lambda_q: C^{s}(\pd \Om)\;\ra\; C^{s-1}(\pd\Om), \quad f\;\ma\; \nu\cdot\nabla u_f|_{\partial\Omega},$$ where $u_f$ is the solution to (\ref{main1}) and $\nu$ is the unit  outer normal vector to  the boundary $\pd \Om.$

  Let $\Gamma_1$ and $ \Gamma_2\subset \partial \Omega$ be two proper open subsets of the boundary. The $W^{2,p}$ well-posedness  guarantees that   the partial boundary Dirichlet-to-Neumann map, denoted by $\Lambda^{\Gamma_1,\Gamma_2}_q$, is well extended to \[\Lambda^{\Gamma_1,\Gamma_2}_q(f):=\Lambda_q(f)|_{\Gamma_2},~ \textrm{ for all } f\in W^{2-\frac1p, p}(\pd \Om) \textrm{ with }\supp (f)\subset \Gamma_1.\]
We are interested in the unique determination of the nonlinear potential $q(x)$ in $\Omega$ from the partial DN map $\Lambda^{\Gamma_1,\Gamma_2}_q$. Moreover, we solve an inverse problem of simultaneously determining the potential and an embedded obstacle or cavity in $\Omega$ from the partial DN map, similar to that in \cite{LLLS2} where the problem was formulated for the semilinear elliptic equation $\Delta u+a(x,u)=0$.

In the second part of the paper, we consider the inverse boundary value problems for two types of nonlinear Maxwell's equations, that arise in potential applications using electromagnetic waves to explore and image the nonlinear properties of the media, using measurements of electromagnetic fields on part of the boundary. The first model of nonlinear electromagnetic medium we consider is for the Kerr-type, described by the boundary value problem
\begin{align}
     \left\{
  \begin{aligned}
  \dt E-i \omega \mu_0 H&=i \omega p(x) |H|^2 H,\\
 \dt H+i \omega \varepsilon_0 E &=-i \omega q(x) |E|^2 E,
  \end{aligned}
           \right.\label{M1}
\end{align}
 in $\Omega\su \R^3,$ with $
    \nt E|_{\partial \Omega}   = f,$
where    $\nu$ is the unit outer normal vector to the boundary $\partial\Omega$. The electrical permittivity $\ep_0$ and the magnetic permeability $\mu_0$ are positive constants and $p(x), q(x)\in C^1(\overline\Omega)$. Here, with some abuse of notations we still denote by $\Omega$ the domain for Maxwell's equations in dimension three instead of arbitrary dimension as considered in the previous Helmholtz equation.

The above   Maxwell's equations of the Kerr-type  are derived from the time-dependent Maxwell's equations
\[\nabla\wedge \mathcal E+\partial_t\mathcal B=0,\quad \nabla\wedge \mathcal H-\partial_t\mathcal D=0,\]
for time harmonic solutions with a fixed frequency $\omega$
\[\mathcal E(x,t)=E(x)e^{-i\omega t}+\overline{E(x)}e^{i\omega t},\quad \mathcal H(x,t)=H(x,t)e^{-i\omega t}+\overline{H(x)}e^{i\omega t},\]
assuming nonlinear constitutive relations of the displacement $\mathcal D$ and induction $\mathcal B$
\[\mathcal D=\varepsilon_0\mathcal E+P_{\textrm{NL}}(\mathcal E),\quad \mathcal B=\mu_0\mathcal H+\mathcal M_{\textrm{NL}}(\mathcal H).\]
Here $\mathcal P_{\textrm{NL}}$ and $\mathcal M_{\textrm{NL}}$ are the nonlinear polarization and magnetization. In a Kerr-type electromagnetic medium, usually that of centrosymmetric structure, the nonlinear polarization is of the form
\[\mathcal P_{\textrm{NL}}(x,\mathcal E(x,t))=\chi_e(x,|E|^2)\mathcal E(x,t),\]
where $\chi_e$ is the scalar susceptibility depending on a memory term such as the time-average of the intensity of $\mathcal E$ \footnote{Note that the time-averages of the intensities of $\mathcal E$ and $\mathcal H$ are $2|E|^2$ and $2|H|^2$ respectively.}. It is common to use $\chi_e(x,|E|^2)=p(x)|E|^2$.
The reader is refereed to \cite{Nie, Stuart} for more details and other examples of electric nonlinear phenomena.
Similarly,  the nonlinearity can be generalized  to magnetization by assuming
\[\mathcal M_{\textrm{NL}}(x,\mathcal H(x,t))=\chi_m(x,|H|^2)\mathcal H(x,t),\quad \chi_m(x,|H|^2)=q(x)|H|^2.\]
Such nonlinear  magnetizations appear in the study of certain metamaterials built by combining an array of wires and split-ring resonators embedded into a Kerr-type dielectric. See \cite{LT} for a numerical implementation of this nonlinear assumption.

The well-posedness of the forward problem of (\ref{M1}) is given  in \cite{AZ1} that except for a discrete set of frequencies, there exists a unique solution $(E, H)\in W^{1,p}_{\textrm{Div}}(\Omega)\times W^{1,p}_{\textrm{Div}}(\Omega)$ for sufficiently small $f$ in a  $L^p$ Sobolev space $ TW^{1-1/p,p}_{\textrm{Div}}(\partial\Omega)$ with $3<p\leq 6$ defined by
\begin{align*}
W^{1,p}_{\textrm{Div}}(\Omega)&=\{u\in W^{1,p}(\Omega;\C^3): \textrm{Div}(\nu\wedge u)\in W^{1-1/p,p}(\partial\Omega)\},\\
TW^{1-1/p,p}_{\textrm{Div}}(\partial\Omega)&=\{f\in W^{1-1/p,p}(\partial\Omega;\C^3): \textrm{Div}(f)\in W^{1-1/p,p}(\partial\Omega)\},
\end{align*}
where $\textrm{Div}$ is the surface divergence operator on $\partial\Omega$. Thus    the admittance map $\Lambda_{p,q}:$ $ TW^{1-1/p,p}_{\textrm{Div}}(\partial\Omega)\rightarrow TW^{1-1/p,p}_{\textrm{Div}}(\partial\Omega)$ is well defined by
\[\Lambda_{p,q}(f):=\nu\wedge H|_{\partial\Omega}.\]
For the corresponding partial data inverse problem, we use the partial admittance map associated to two open non-empty open sets $\Gamma_1$ and $\Gamma_2$ of $\partial\Omega$ as
\[\Lambda^{\Gamma_1,\Gamma_2}_{p,q}(f):=\Lambda_{p,q}(f)|_{\Gamma_2}=\nu\wedge H|_{\Gamma_2}\]
for all $f\in C_c^s(\Gamma_1;\C^3)\cap TW^{1-1/p,p}_{\textrm{Div}}(\partial\Omega)$. We consider the inverse problem in determining the nonlinear coefficients $p$ and $q$ from this partial admittance map.

Another important nonlinear electromagnetic behavior of the materials is second harmonic generation (SHG). Applications of the nonlinear optical phenomena includes obtaining coherent radiation at a wavelength shorter than that of the incident laser, through the frequency doubling effect of SHG. Moreover, in the {\em second-harmonic imaging microscopy (SHIM)}, a second-harmonic microscope obtains contrasts from variations in {\em a specimen's ability to generate second-harmonic light} from the incident light while a conventional optical microscope obtains its contrast by detecting variations in optical density, path length, or refractive index of the specimen.
The SHIM is also exploited in imaging flux residues (see the work in Chen Lab at the University of Michigan). Although nonlinear optical effects are in general very weak, the significant enhancement of SHG was shown using diffraction gratings or periodic structures.

We consider an inverse boundary value problem for Maxwell's equations
\begin{align*}
\left\{
\begin{array}{ll}
\dt E^{\om}-i \omega \mu_0 H^{\om}&=0\\
 \dt H^{\om}+i \omega \varepsilon_0 E^{\om} &=-i \omega \chi^{(2)} \overline{E^{\om}}\cdot E^{2\om}\\
  \dt E^{2\om}-i2 \omega \mu_0 H^{2\om}&=0 \\
 \dt H^{2\om}+i 2\omega \varepsilon_0 E^{2\om} &=-i   2\omega \chi^{(2)}    E^{\om} \cdot E^{ \om}\\
 \nt E^{ \om}|_{\pd \Om} &=f^{ \om}\\
  \nt E^{ 2\om}|_{\pd \Om} &=f^{ 2\om}
 \end{array}
 \right.
\end{align*}
  where  $f^{k\om}\in C^{s}_c(\pd \Om)\cap TW_{\textrm{Div}}^{1-1/p,p}(\partial\Omega)$ and  $ \supp f^{k\om}\su \Gamma_1 $ for $k=1,2$. This models the phenomenon when a beam with time-harmonic electric field
$E(t, x)=E(x)e^{-i\omega t} + \textrm{c.c.}$
is incident upon an SHG medium (e.g., a noncentrosymmetric crystal), new waves are generated at zero frequency and at frequency $2\omega$, and assuming that the susceptibility parameter is isotropic $\chi^{(2)}=(\chi^{(2)}_l)_{l=1}^3$.
The corresponding  output of the admittance map
$$\Lambda_{\chi^{(2)}}^{\Gamma_1,\Gamma_2}(f^{k\om})=\nt H^{k\om}|_{\Gamma_2}$$ is measured on $\Ga_2$ with $\Gamma_1\cap \Gamma_2\not=\emptyset$. $\chi^{(2)}$ is the second order susceptibility parameter. Then we are solving the problem of determining the susceptibility $\chi^{(2)}$. \\

The type of inverse boundary value problems was first formulated by Calder\'{o}n in \cite{Ca80} for a proposed imaging method, known as the electrical impedance tomography, in which one aims to
  determine the electrical conductivity, modeled by the function $\sigma(x)$ for $x$ in a body $\Omega$, from the boundary measurements of the electric voltage and current, formulated using the DN map $\Lambda_\sigma: u|_{\partial\Omega}\mapsto \nu\cdot\sigma\nabla u|_{\partial\Omega}$ for the linear elliptic equation $\nabla\cdot(\sigma\nabla u)=0$ in $\Omega$. The first global uniqueness result was proved by Sylvester and Uhlmann    \cite{SU87}   for $C^2$-conductivities $\sigma(x)$
    in dimensions $n\geq 3$ by solving the problem of determining an electric potential $q(x)$ in the Schr\"odinger equation $(-\Delta+q)u=0$ from the boundary DN map. 
      Later the  regularity condition of conductivities   was relaxed in \cite{BT} and \cite{PPU}.
 In  two dimensions,  Nachman in  \cite{Na96} and Astala-P\"{a}iv\"{a}rinta in  \cite{AP} proved the global uniqueness  result   for $C^2$-conductivities and    for $L^{\infty}$-conductivities, respectively. Bukhgeim in \cite{Bu} obtained the uniqueness result for $L^p$-potentials of Schr\"{o}dinger equations  
  from Cauchy data in two dimenions.
 Many further results for global uniqueness of inverse boundary value problems   are available in the literature. We refer readers to
  \cite{ALP,S90,SU91,Sun1} for anisotropic conductivities in two and higher dimensions, \cite{NaSU,Sa} for magnetic Schr\"{o}dinger operators and  \cite{OPS93,OS96,KSU} for Maxwell's equations, and \cite{DKSU2,DKuLS,DKuLS2,FO} in Riemannian  geometries.  We also refer readers to   the surveys  \cite{Uh09,Uh14} and the references therein.


The study of the inverse problem for nonlinear elliptic equations goes back to the 1990's. The method introduced by Isakov in \cite{Isa93} for the parabolic equations, shows that the linearization of the nonlinear DN map indeed gives full information about the DN map of the corresponding linear equation, and thus the
  uniqueness results of inverse problems for linear equations could be applicable. The method was then generalized to recover nonlinear coefficients for a class of   semilinear equations  \cite{Isa94,IN95,Sun10} in two and higher dimensions and   for a parabolic systems of semilinear equations in \cite{Isa01}. Also see \cite{IsaCPDE} for an application of the linearized method to quasilinear equations to recover  nonlinear coefficients from partial data.
Furthermore, a second order linearization was applied to a class of   quasilinear equations  for the unique determination of  quadratic terms in \cite{KaNa}  and determination of anisotropic conductivities in \cite{Sun,SunU}.  Also see   \cite{AZ1,AZ2,FO,KU1,KU2,LLLS1,LLLS2} for an application of higher order linearization methods  to nonlinear elliptic equations  for the recovery of    power-type nonlinear terms.
More details on inverse problems for nonlinear equations can be found in the surveys \cite{Uh09,Uh14}.


In the last few years, the higher order linearization method has become a powerful tool in dealing with inverse problems for nonlinear hyperbolic equations including     wave equations in  \cite{FO2,KLU,LUW1,HUW,WZ}  and einstein equations in \cite{LUW2,KLOU}. In those works,  the nonlinearity is proven to
be helpful in determining the information of the coefficients in the operators, in some situations combined with microlocal analysis of the newly generated singularities in order to
solve the inverse problems whose corresponding version for linear equations are otherwise still open. \\


Our first result for the nonlinear Helmholtz Schr\"odinger equation is given below.
\begin{theo}\label{uq1}
Let $q_1$ and $q_2\in L^{\infty}(\bar{\Omega}),$ and the Dirichlet-Neumann map $\Lambda_{q_j}$ satisfy that
\begin{align*}
  \Lambda_{q_1}(f)|_{\Gamma_2}=\Lambda_{q_2}(f)|_{\Gamma_2},~~\forall f\in C_c^s(\Gamma_1),
\end{align*}
with $\|f\|_{C_c^s(\Gamma)}<\delta$, $s>1$ with $s\not\in\Z$, where $\delta$ is sufficiently small. Then $q_1=q_2$ in $\Omega.$
 \end{theo}

\begin{rem}
Following Theorem \ref{aw1} in appendix, the well-posedness of nonlinear Helmholtz equation  holds for small boundary value $ f\in W^{2-\frac1p,p}(\pd \Om) .$ Thus the  Dirichlet-to-Neumann map is well defined and uniquely determines $q(x)\in L^{\infty}(\Om)$. The well-posedness is also given in \cite{LLLS1} in H\"older spaces. 
\end{rem}

We remark here that our result for the nonlinear Helmholtz type Schr\"odinger equation \eqref{main1} is not contained in the cases discussed in \cite{LLLS2}, where the linearized equation is of Laplace type. Due to the non-positivity of operator $-\Delta-k^2$, one needs modified density arguments discussed below. On the other hand, the equation \eqref{main1} can be viewed as a fixed frequency time-harmonic equation for the wave operator $\partial_t^2-\Delta$. Therefore, our result sheds some light on the inverse boundary problem for wave equations. In \cite{DKuLS}, the elliptic DN map of Calder\'{o}n problem of Schr\"{o}dinger equation in an infinite cylinder
 is reduced to the hyperbolic DN map of a  wave equation. Based on the unique continuation property for  the reduced wave equation, the boundary control method is then applied to recover the uniqueness. Our paper presents a direct proof of partial data inverse problem for Helmholtz Schr\"{o}dinger equation. We also remark a progress in
   \cite{OSSU}, which reduces Calder\'{o}n problems to the injectivity of some geodesic X-ray transforms and   and provides a unified approach for inverse boundary problems for Laplace type, transport and wave equations.

  In general,
 the higher order linearization method  gives rise to an integral identity which involves an  integral of some product of solutions to the linearized equations. The problem is then reduced to an integral geometry problem after plugging in proper linear solutions. When the linearized equation is the Laplace equation and the boundary DN-map is given as partial data $\Lambda^{\Gamma_1,\Gamma_2}$, the solutions are usually chosen to be harmonic exponentials that vanish on part of the boundary. The density result of the product of two such exponentials was shown in \cite{DKSU} as an enlightening step  of the study of   Calder\'on problem for Schr\"odinger equation $-\Delta u+qu=0$.   Their techniques are further exploited in \cite{KU1} to show the density result of the set of the product of two gradients of some harmonic functions with some  vanishing boundary values. Such density results were then directly or indirectly used in proving the uniqueness of   inverse boundary problems for nonlinear elliptic equations in \cite{ KU1, KU2, LZ,LLLS1, LLLS2}. As a comparison, for the partial data problems of linear equations,  two approaches are mainly proposed: one is a reflection argument in \cite{Isa07,COS} by assuming that the inaccessible boundary  is either on a plane or on a sphere; the other is to establish the Carleman estimate as in \cite{KSU,DKSU2} where the partial  data on the boundary are closely related with the constructed limiting Carleman weight. A combination of  these two approaches is exploited in \cite{KS} to extend the uniqueness result to manifolds. These methods aim at constructing new complex geometrical optics (CGO) solutions similar to those in \cite{SU87} for the linear equations and usually requires some extra assumptions on $\Gamma_1$ and $\Gamma_2$.

In  our paper, following the techniques in \cite{DKSU}   we first prove the following density result in order to obtain the global uniqueness of the coefficient for a nonlinear Helmholtz Schr\"{o}dinger equation after higher order linearization. 

\begin{theo}\label{ds1}
   Let $\Omega$ be a bounded open set in $\R^n$,   $n\geq 2$ with smooth boundary $\partial \Omega$, and $\Gamma\subset \partial \Omega$ be a proper nonempty open subset of the boundary. Let $q(x)\in L^{\infty}(\Om).$ \\ Suppose the cancellation equality
 \beq\label{ce1}
  \int_{\Omega} q(x)v_1(x)v_2(x)dx=0
\eeq
  holds for any smooth solutions $v_1(x), v_2(x)\in C^{\infty}(\bar{\Omega})$ to the equation
  \begin{align}\label{ce2}
\left\{\begin{array}{ll}
  -\Delta v(x)-k^2 v(x) =0&~\text{in~}\Omega,\\
     v(x)=0&\text{on~} \partial \Omega\setminus\Gamma.
   \end{array}
   \right.
\end{align}
  Then   $q(x)$ vanishes in $ \Omega$.
\end{theo}

As another application of the density result,  our paper also yields the simultaneous recovery of unknown cavity  or boundary and coefficients  for  the nonlinear Helmholtz Schr\"{o}dinger operator $-\Delta u+k^2 u +q(x)u^2$ following the discussion in \cite{LLLS2}. In comparison, the solutions to the linearized equation  in our case are no longer harmonic functions, and   the  tools available  are restricted to  the unique continuation principle.\\







Our main results for the nonlinear Maxwell's equations are given below.

\begin{theo}[Kerr-type]\label{KL1}
  Let $p_1,p_2, q_1$ and $q_2\in L^{\infty}(\bar{\Omega})$ and $\Gamma_1\cap\Gamma_2\neq\emptyset$. Suppose
\begin{align*}
  \Lambda_{p_1,q_1}(f)|_{\Gamma_2}=\Lambda_{p_2,q_2}(f)|_{\Gamma_2},~~\forall f\in C_c^s(\Gamma_1;\C^3)\cap TW_{\textrm{Div}}^{1-1/p,p}(\partial\Omega),
\end{align*}
with $\|f\|_{C_c^s(\Gamma)}<\delta$, $s>1$ with $s\not\in\Z$, where $\delta$ is sufficiently small. Then $p_1=p_2$  and $q_1=q_2$ in $\Omega.$
\end{theo}

We also give the identifiability result of the  partial data   inverse problem of a class of  nonlinear time-harmonic Maxwell's system with the second harmonic generation
 on $\Om$

\begin{theo}[Second Harmonic Generation]\label{HL1}
  Let $\chi^{(2),1} $ and $\chi^{(2),2}\in ( L^{\infty}(\bar{\Omega}))^3,$ the admittance map $\Lambda_{\chi^{(2)}}$ satisfy that
\begin{align*}
  \Lambda_{\chi^{(2),1} }(f^{k\om})|_{\Gamma_2}=\Lambda_{\chi^{(2),2}}(f^{k\om})|_{\Gamma_2},~~\forall f^{k\om}\in C_c^s(\Gamma_1)\cap TW_{\textrm{Div}}^{1-1/p,p}(\partial\Omega),
\end{align*}
with $\|f^{k\om}\|_{C_c^s(\Gamma)}<\delta$, $s>1$ with $s\not\in\Z$, where $\delta$ is sufficiently small for $k=1,2$. Then $\chi^{(2),1}=\chi^{(2),2}$  in $\Omega.$
\end{theo}

  For partial data inverse problems  of Maxwell's equations, the literature relies on geometrical assumptions   on  the inaccessible boundary $\Ga_2^c$ to be either on a plane or on a sphere  \cite{COS} or on the admissible  manifold  \cite{COST} with     $\Ga_1$ to be the global boundary. Our result applies the high order linearization method and   $\Ga_1$ and $\Ga_2$ could be arbitrary.

We first extend the  density result in \cite{DKSU} which states that the set of a product of harmonic functions with boundary value supported on part of the boundary of the domain is dense in $L^2(\Om)$ to the set of a product of solutions to linear Helmholtz equations with boundary value supported on part of the boundary.  Different from calculating symbols of the nonlinear interaction, we adopt integral identity  and follow the microlocal analysis as in \cite{DKSU} to make the proof simple.
Then we apply the result to recover the coefficients of a class of nonlinear Helmholtz Schr\"{o}dinger and Maxwell's equations with partial data by virtue of the higher order linearization method. We also give some simultaneous recovery results of coefficients and unknown cavity  as applications. \\

The paper is organized as follows. In Section 2 we give the proofs of  Theorem \ref{uq1} and Theorem \ref{ds1}. Additionally, we provide some applications of simultaneous recovery of  nonlinear  Helmholtz Schr\"{o}dinger equations in Section 3. Section 4 is devoted to the identifiability results of Maxwell equations, i.e.,  Theorem \ref{KL1} and Theorem \ref{HL1}.

\section{Nonlinear Helmholtz Schr\"odinger equations }
\subsection{Proof of Theorem \ref{uq1} }\label{yy1}
Following \cite{LLLS1,LLLS2}, we adopt the higher order linearization of the DN-map to deal with the nonlinearity. Since the application of  the high order linearization method here is standard, we provide key procedures   and refer readers to \cite{LLLS1,LLLS2} for details.

\Proof
Let $f=\epsilon_1f_1+\epsilon_2f_2$, where $\epsilon_1, \epsilon_2>0$, and  $f_1$ and $f_2\in C_c^s(\Gamma_1).$ Consider the parameterized boundary value problem
\begin{align}\label{sec21}
\left\{\begin{array}{lll}
  -\Delta u_j-k^2 u_j+q_j(x)u_j^2 =0 &~\text{in~}\Omega,\\
   u_j  =\epsilon_1f_1+\epsilon_2f_2  &~\text{on~}\pd\Omega,
   \end{array}
   \right.
\end{align}
with $j=1,2$.
It follows from the well-posedness of the nonlinear Helmholtz  equation given in Appendix \ref{Appendix:well-posedness} that   when $k$ is not an eigenvalue of the Laplace operator,  there exist unique solutions $ u_j:=u_j(x;\epsilon_1, \epsilon_2) $ to (\ref{main1}) with $j=1,2$, provided $\epsilon_1, \epsilon_2$ are sufficiently  small. Then the assumption $\Lambda_{q_1}(\epsilon_1f_1+\epsilon_2f_2  )|_{\Gamma_2}=\Lambda_{q_2}(\epsilon_1f_1+\epsilon_2f_2)|_{\Gamma_2}$  gives
 \beq\label{dn1}
 \partial_{\nu}u_1|_{\Gamma_2}=\partial_{\nu}u_2|_{\Gamma_2}.
 \eeq


Denote $v_j^{(l)}:= \partial_{\epsilon_l} u_j\big|_{\varepsilon_l=0}$ for $j=1,2$ and $l=1,2.$ The  differentiation of (\ref{sec21}) w.r.t. $\epsilon_l$ at $\e_1=\e_2=0  $ and the initial value $u_j(x;0,0)=0$ yield
\begin{align*}
\left\{\begin{array}{ll}
  -\Delta v_j^{(l)}-k^2 v_j^{(l)}   =0 &~\text{in~}\Omega,\\
   v_j^{(l)}  = f_l &~\text{on~}\pd\Omega,
   \end{array}
   \right.
\end{align*}
with $j=1,2$ and $l=1,2.$
 $v_1^{(l)}=v_2^{(l)}$ follows from the uniqueness of  Dirichlet boundary value problem of Helmholtz equation    for each $l$.
From now on,
denote by
\begin{align*}
 v^{(l)}:= v_1^{(l)}=v_2^{(l)},~\text{with~}l=1,2;
\end{align*}
thus $ v^{(l)}$ is a solution to
\begin{align}\label{vjl1}
\left\{\begin{array}{ll}
  -\Delta v ^{(l)}-k^2 v ^{(l)}    =0 &~\text{in~}\Omega,\\
   v ^{(l)} = f_l &~\text{on~}\pd\Omega.
   \end{array}
   \right.
\end{align}
   Denote $w_j := \partial^2_{ \epsilon_1   \epsilon_2}u_j \big|_{\varepsilon_1=\e_2=0}$, with $j=1,2 $.
Then the second order linearization of (\ref{sec21}) yields
\begin{align*}
\left\{\begin{array}{ll}
  -\Delta  w_j -k^2 w_j  +2 q_j(x) v^{(1)} v^{(2)} =0& ~\text{in~}\Omega,\\
    w_j  = 0 & ~\text{on~} \pd \Om.
   \end{array}
   \right.
\end{align*}
Moreover, from (\ref{dn1}) one has
\beq\label{w1}
\partial_{\nu}w_1 |_{\Gamma_2}=\partial_{\nu}w_2 |_{\Gamma_2}.
\eeq

\vspace{5mm}
\noindent

In order to establish an integral identity  we choose some proper solution $v^{(0)}$ to compensate the unknown information on $\Gamma_2$ .
Let $v^{(0)}$ be a solution to
\begin{align}\label{v0}
\left\{\begin{array}{ll}
  -\Delta v^{(0)}-k^2 v^{(0)}   =0&~\text{in~}\Omega,\\
    v^{(0)} = f_0   & \text{on~} \Gamma_2,
   \end{array}
   \right.
\end{align}
where $f_0\in C_c^s(\Gamma_2)$ and recall that  $\Gamma_2 \su \pd \Om$ is an open nonempty subset.
Thus,
\begin{align*}
   & 2 \int_{\Omega}  (q_1(x)-q_2(x))v^{(1)}v^{(2)} v^{(0)}dx\\
    =&2  \int_{\Omega}    \Delta (w_1-w_2)v^{(0)} +  k^2  (w_1-w_2)   v^{(0)} dx\\
  =& 2 \int_{\Omega} \Delta (w_1-w_2)v^{(0)} -\Delta v^{(0)} (w_1-w_2) dx\\
   = & 2 \left(\int_{\partial \Omega} \frac{\partial}{\partial \nu} v^{(0)}(w_1-w_2)dS- \int_{\partial \Omega} \pd(w_1-w_2) v^{(0)} dS\right)
\end{align*}
vanishes in view of  (\ref{w1}) and (\ref{v0}).
Hence   the integral identity
\beq\label{ii2}
\int_{\Omega}  (q_1(x)-q_2(x))v^{(1)}v^{(2)} v^{(0)}dx=0,
\eeq
holds, where  $v^{(1)}, v^{(2)}$ are solutions to the linearized equation  (\ref{vjl1}) and  $ v^{(0)}$ satisfies the     boundary value problem   (\ref{v0}).

Apply Theorem \ref{ds1} and one has
$$
 (q_1(x)-q_2(x))v^{(0)}(x)=0
$$
identically in $\Omega$, where $ v^{(0)}$ is a solution to the     boundary value problem   (\ref{v0}). By constructing a solution  $v^{(0)}$ satisfying  $v^{(0)}(\tilde{x}_0)\not=0 $ for any fixed $\tilde{x}_0\in \Omega$, we are able to show $(q_1-q_2)(\tilde{x}_0)=0$.
In fact, apply  Proposition (\ref{pen2}) with $L=-\Delta-k^2$ and consider  $u(x)=e^{-i x\cdot \xi}$ with $\xi\cdot \xi=k^2$ such that $Lu=0$.
  Moreover,    $\xi$ is chosen such that $u(\tilde{x}_0 )\not=0$.
   Then  Proposition (\ref{pen2}) yields  a solution $v^{(0)}(x)$   close to $u(x)$ and thus $ v^{(0)}(\tilde{x}_0)\not=0$ by the continuity.
   Since $\tilde{x}_0$ is arbitrary in $\Omega$,  $q_1=q_2 $ in $\Omega$. This completes the proof of Theorem   \ref{uq1}.
\hfill $\Box$

\subsection{Proof of Theorem \ref{ds1} }
The proof of Theorem \ref{ds1} consists of  three parts. First we   establish the new setting of the local problem  and prove  local results of   Theorem \ref{ds1} in section \ref{sh1}.
Then section 2.2.2 contributes to the    extension of   local results to global results is given   in section.

\subsubsection{Local results  }\label{sh1}

\begin{prop}\label{loc1}
  Let $\Omega$ be a bounded open set in $\R^n$, with $n\geq 2$ with smooth boundary $\partial \Omega$,  $x_0\in \partial \Omega$   a convex point on the boundary, and
   $\Gamma$ be  a small neighborhood of $x_0$ on the boundary. Suppose that the cancellation equality (\ref{ce1}) holds for any
 smooth  solutions $v_1(x)$ and $ v_2(x)\in C^{\infty}(\bar{\Omega})$ to the equation
(\ref{ce2}). Then there exists $\delta>0$ such that $q(x)$ vanishes on $B(x_0,\delta)\cap \Omega$.
\end{prop}

\Proof
Choose a point $a\in \R^n\setminus \bar{\Omega}$ along the direction of the inward normal vector at $x_0\in \partial \Omega$. As $x_0$ is convex  and  $\Omega$ is bounded, suppose that there exists a ball $B(a,r)$ containing $\Omega$ and $\partial B(a,r)\cap \bar{\Omega}=\{x_0\}$.
Then   a transform    $  \varphi(x)$  consisting of   translation and  rotation transforms   such that $x_0$ is mapped  to the origin and $a$ is mapped to $-e_1=(-1,0,\ldots, 0)$.
Indeed,
  $ \varphi $ can be defined as
 $  x\in  \Omega \mapsto    y = \frac{1}{r}\cdot  O(x-a) -e_1,$ where $O$ is an  orthogonal matrix such that $O(x_0-a)=|x_0-a|\cdot e_1=r e_1.$
 Thus $\tilde{\Omega}=\varphi (\Omega)$ is contained in $B(-e_1,1).$

Note that the Laplacian $\Delta$ is invariant under $\varphi$ as a combintion of translation and rotation transform.
Then, denote  $\tilde{q}=q \circ \varphi^{-1} ,$
  $\tilde{u}_1(y)= u_1\circ \varphi^{-1}$ and $\tilde{u}_2(y)=u_2  \circ \varphi^{-1}$. The   cancellation (\ref{ce1})  is changed to
$
 \int_{\tilde{ \Omega }} \tilde{q}   \tilde{u}_1    \tilde{u}_2     dy=0,
$
where $\tilde{u}_1$ and $\tilde{u}_2$ satisfy   $  -\Delta u-k^2r^2 u=0$ in $\tilde{ \Omega } $  and vanish on the boundary $\tilde{\Omega}\setminus \tilde{\Gamma}$.

If $\tilde{q}$ vanishes near $y_0=0$, then $q$ also vanishes near $x_0$. Therefore,  we start with a new cancellation equality. With some abuse of notations we omit the tilde   sign     and still use $x$-coordinate. In our new setting,    denote
  $x_0=0.$      Suppose  that
  $  |x+e_1|\leq 1$  for $x\in \Omega,$
and
$
 \partial\Omega \setminus \Gamma \subset \{x\in \partial \Omega |~ x_1\leq-2  c \}.
$
The following cancellation holds
\begin{align*}
 \int_{  \Omega  } q  u_1    u_2     dx=0,
\end{align*}
 for
any smooth  solutions  $u_1$ and $u_2$ to the equation
    \begin{align*}
\left\{\begin{array}{ll}
    -\Delta u-k^2 r^2 u=0&~\text{in~}\Omega,\\
    u  = 0 & \text{on~} \partial \Omega\setminus \Gamma.
   \end{array}
   \right.
\end{align*}

The next proof  of local results  follows the literature \cite{DKSU}. Thus to shorten the length we only provide some critical procedures and estimates here. For more details of the calculation we refer to \cite{DKSU}.

We first construct $u_1$ and $u_2$  by adding a correction term  $w_1(x,\zeta)$ and $w_2(x,\eta)$   to eliminate  the value of CGO solutions  of the form  $e^{-\frac{i}{h}x\cdot \zeta } $  and $e^{-\frac{i}{h}x\cdot \eta }$  on the boundary $ \Omega \setminus \Gamma $, respectively, where $\zeta$ and $\eta$ are   proper complex vectors.

Note that $p(\xi)=\xi^2$ is the principal symbol of the Laplacian on $\R^n$ and consider its complexification  on $\C^n$. Let  $\zeta$ and $\eta$ belong to the set
$$
p^{-1}(krh)=\{\zeta\in \C^n|~\zeta\cdot \zeta=(k r h)^2\}.
$$
 Choose $\zeta_0=(i,\sqrt{(k r h)^2+1},0,\ldots,0)$  and  $\eta_0=(i,-\sqrt{(k r h)^2+1},0,\ldots,0).$  Then $\zeta_0$ and $\eta_0$ belong to $p^{-1}(krh).$
Since the differential of the map
$s:p^{-1}(krh)\times p^{-1}(krh)\ra \C^n, $ $(\zeta,\eta)\mapsto \zeta+\eta$ at $(\zeta_0,\eta_0)$ is surjective, any complex  vector $z\in \C^n$ with $|z-2ai e_1|<2 a \varepsilon$ can be decomposed as a sum of the form
\beq\label{dec1}
z=\zeta+\eta,~~~\zeta,\eta\in p^{-1}(krh),~|\zeta-a\zeta_0|<Ca\varepsilon,|\eta-a\eta_0|<Ca\varepsilon.
\eeq

In the sequel, we use $C$ to denote different constants independent of $h$ in the inequalities for convenience.

Take a cutoff function $0\leq \chi\leq 1\in C^{\infty}_c(\R^n)$   which equals 1 on $x_1\leq -2c$ containing  $\partial\Omega\setminus \Gamma$
with  compact  support  in $x_1\leq -c$. There exists a  solution $w_1(x,\zeta)$ to the Dirichlet problem
 \begin{align}\label{cor1}
\left\{\begin{array}{ll}
  -\Delta  w -(kr)^2 w =0&~\text{in~}\Omega,\\
    w    =  - \chi(x)e^{-\frac{i}{h}x\cdot \zeta }&~\text{on~}\pd\Omega.
    \end{array}
   \right.
\end{align}
Thus, fix $k$ outside a discrete set $\Sigma_r$, one has
$$
\| w_1\|_{H^1 (\Omega)}  \leq    C
\| \chi(x)e^{-\frac{i}{h}x\cdot \zeta }\|_{H^{\frac32} (  \partial \Omega ) }
$$
and
\begin{align}\label{crz1}
\| w_1\|_{H^1 (\Omega)}  \leq
 C(1+h^{-1}|\zeta|)^{\frac32}e^{-\frac ch \im \zeta_1 }e^{ \frac 1h |\im \zeta^{\prime}| },~~\text{if~} \im  \zeta_1\geq 0.
\end{align}
Let $u_1(x,\zeta)=e^{-\frac{i}{h}x\cdot \zeta } +w_1(x,\zeta)$, with $\zeta\in p^{-1}(krh)$. Then $u_1$ is a desired solution and   vanishes on $\pd \Om\setminus \Gamma$. Similarly, we construct $u_2(x,\eta)=e^{-\frac{i}{h}x\cdot \eta } +w_2(x,\eta)$, where $\eta\in p^{-1}(krh)$ and  $w_2(x,\eta)$ is a correction term satisfying a similar estimate to \eqref{crz1} with $\zeta$ replaced by $\eta$.
Plug in  $u_1(x,\zeta)$ and   $u_2(x,\eta)$ into the cancellation and combine the   estimates of $\| w_1\|_{H^1 (\Omega)} $ and $\| w_2\|_{H^1 (\Omega)} $. When $ \im  \zeta_1\geq 0$ and $ \im  \eta_1\geq 0$, the following estimate
\begin{align*}
  \left|\int_{  \Omega  } q(x)  e^{-\frac ih x \cdot (\zeta+\eta) } dx\right|
\leq     C \| q\|_{L^{\infty}(\Omega)}& (1+h^{-1}|\zeta|)^{\frac32}(1+h^{-1}|\eta|)^{\frac32}\\
& e^{-\frac ch \min\{\im \zeta_1, \im \eta_1\} }
e^{ \frac 1h( |\im \zeta^{\prime}|+ |\im \eta^{\prime}|) },
\end{align*}
 holds.
In particular, when
$|\zeta-a\zeta_0|<Ca\varepsilon $ and $|\eta-a\eta_0|<Ca\varepsilon$ with $\varepsilon\leq \frac{1}{2C}$,
\begin{align*}
  \left|\int_{  \Omega  } q(x)  e^{-\frac ih x \cdot (\zeta+\eta) } dx\right|
\leq     C h^{-3} \| q\|_{L^{\infty}(\Omega)}
e^{-\frac {ca}{2h}   }
e^{ \frac {2Ca \varepsilon}{h} }.
\end{align*}
   Therefore, for all $z\in \C^n$ with $|z-2  a i e_1|<2 a \varepsilon$, it follows from the decomposition (\ref{dec1}) that
\begin{align*}\label{fbi1}
   \left|\int_{  \Omega  } q(x)  e^{-\frac ih x \cdot z }  dx\right|
\leq   C h^{-3} \| q\|_{L^{\infty}(\Omega)}
e^{-\frac {ca}{2h}   }
e^{ \frac {2Ca \varepsilon}{h} }.
\end{align*}
This is a similar estimate as in \cite{DKSU}, which guarantees
  an exponential decay for  the F.B.I. transform of $f $.
Finally, an  application of  the water melon approach as in \cite{DKSU} yields
%
%
 $$
q(x)=0,~~\forall x\in \Omega,~~-\delta\leq x_1\leq 0,
 $$
provided $\delta$   small enough. This completes the proof of Proposition \ref{loc1}. \hfill $\Box$

\subsubsection{Extend to global results }
In this section, we extend the local result (Proposition \ref{loc1}) globally to the whole $\Om$. The following lemma is needed.

\begin{lem} \label{loc3}
Let $\Omega_1\subset\Omega_2$ with smooth boundary. Define $G_{\Omega_2}(x,y)$ to be the Green function of the following system
  $$
  \left\{
  \begin{array}{ll}
  (-\Delta_y-k^2)G_{\Omega_2}(x,y) =\delta(x-y)  &~\text{in ~}  \Omega_2,\\
  G_{\Omega_2}(x,y) =0 &~\text{on ~} \pd  \Omega_2.
  \end{array}
  \right.
  $$
  The set
  $$
  \mathcal{R}=\left\{
  \int_{\Omega_2}  G_{\Omega_2}(x,y)a(y)dy,~\text{for some~}a(y)\in C^{\infty}(\bar{\Omega}_2),~\supp a\subset \bar{\Omega}_2\setminus \Omega_1,x\in \Omega_2
  \right\}
  $$
  is dense for $L^2(\Omega_1)$ topology in the  space  $\mathcal{S}=\{u(x)\in C^{\infty}(\bar{\Omega}_1):
    (-\Delta_x-k^2)u=0,u|_{\partial \Omega_1\cap\partial \Omega_2 }=0.\}$
\end{lem}
\begin{Proof}
    It suffices to show that  $\lag v, \mathcal{S} \rag=0$ for any $v(x)\in L^2(\Omega_1) $ such that $\lag v, \mathcal{R} \rag=0$.
  Following   Fubini Theorem and the assumption, the equality
  $$
\int_{\Omega_1} v(x) \int_{\Omega_1}   G_{\Omega_2}(x,y)a(y)dy dx= \int_{\Omega_2}\int_{\Omega_1}  v(x) G_{\Omega_2}(x,y)dx~ a(y)dy =0,
  $$
holds  for any $a(y)\in C^{\infty}_c(\bar{\Omega}_2)$.
The assumption  $\supp a\subset \bar{\Omega}_2\setminus \Omega_1$ yields
 \beq\label{loc2}
 \int_{\Omega_1}  v(x) G_{\Omega_2}(x,y)dx=0,~~\text{for~} y\in \bar{\Omega}_2\setminus \Omega_1.
 \eeq
Denote
$
\omega(y)=   \int_{\Omega_1}  v(x) G_{\Omega_2}(x,y)dx.
$
Then
 $  \omega(y)=0$ for $y\in \bar{\Omega}_2\setminus \Omega_1. $
Moreover, $\omega\in H_0^1(\Omega_2) $ in view of  $G_{\Omega_2}(x,y)|_{y\in \partial \Omega_2}=0$.

Note that
 $ (-\Delta_y-k^2)\omega(y)=v(y) $ for $y\in \Omega_2$.
 Thus, for any $u\in \mathcal{S},$
 \begin{align*}
 \int_{\Omega_1} uvdx&=  \int_{\Omega_1} u(-\Delta_x-k^2)\omega(x)dx=  \int_{\Omega_1}  (\Delta u \omega -u \Delta \omega)dx \\
 &=\int_{\partial\Omega_1\cap \partial\Omega_2}  ( u \pn\omega - \omega \pn u)dS
 +\int_{\partial\Omega_1\setminus \partial\Omega_2}  ( u \pn\omega - \omega \pn u)dS.
 \end{align*}
 The  right hand side term vanishes by combining $u|_{\partial\Omega_1\cap \partial\Omega_2}=0$, $\omega \in H_0^1(\Omega_2)$   and  $\omega=0 $ in $\bar{\Omega}_2\setminus \Omega_1$.
 Hence $ \int_{\Omega_1} uvdx=0$ for any $u\in \mathcal{S}.$
$\hfill\Box$

\end{Proof}

\vspace{5mm}
We are ready    to conclude the proof of Theorem \ref{ds1}. The approach goes through as in \cite{DKSU}.

\noindent
\textit{Proof of Theorem \ref{ds1}~~}
Fix a point $x_1\in \Omega$ and let $\theta:[0,1]\ra \bar{\Omega}$ be a $C^1$ curve joining $x_0\in \Gamma$ to $x_1$ such that $\theta(0)=x_0$  and $\theta^{\prime}(0)$ is the interior normal to $\partial \Omega$ at $x_0$ and $\theta(t)\in \Omega$ for all $t\in (0,1].$
Let
 $
\Theta_{\e}(t)=\{x\in \bar{\Omega}|~d(x,\theta([0,t]))\leq \e\}
 $
be  a closed neighborhood of the curve $\theta(t),$ $t\in [0,1]$, and set
$
I=\{t\in [0,1]|~f=0 \text{~a.e. on~} \Theta_{\e}(t)\cap \Omega \}.
$
 $I$ is a closed set. From Proposition \ref{loc1}, $I$ is non-empty provided that $\e$ is small enough. We would like to  show that $I$ is also open. Thus by connectivity $I=[0,1].$ Then $x_1\not\in\supp f $. Since $x_1$ is arbitrary, $f=0 $ in $\Omega$. Hence, it suffices to show that $I$ is open.

 Let  $t\in I$ and $\e$ small enough. Then we can assume that
 $\partial \Theta_{\e} (t)\cap \partial \Omega\subset \Gamma,
 $
 and $\Omega\setminus \Theta_{\e} (t)$ can be smoothed out into an open subset $\Omega_1\subset \Omega$ with smooth boundary such that
 $
 \Omega\setminus \Theta_{\e} (t)\subset \Omega_1,~~\partial \Omega\setminus \Gamma\subset \partial \Omega\cap \partial \Omega_1.
 $
 Furthermore,
we augment $\Omega$ by smoothing out the set $\Omega\cup B(x_0,\e )$ into an open set $\Omega_2$ with smooth boundary such that
$
\partial \Omega\setminus \Gamma\subset \partial \Omega\cap \partial \Omega_1\subset  \partial \Omega\cap \partial \Omega_2.
$

Denote by $G_{\Omega_2}$ be the fundamental solution   to the  system
$$
  \left\{
  \begin{array}{ll}
  (-\Delta_y-k^2)G_{\Omega_2}(x,y) =\delta(x-y)& ~\text{in ~}  \Omega_2,\\
  G_{\Omega_2}(x,y)  =0&~\text{on ~} \pd \Omega_2,
  \end{array}
  \right.
  $$
and let
$$
H(x,t)=\int_{\Omega_1}  f(y) G_{\Omega_2}(x,y) G_{\Omega_2}(t,y) dy,~~\text{for~} t,x\in \Omega_2\setminus \bar{\Omega}_1.
$$
Then $H$ satisfies  $-\Delta_x H-k^2 H=0$   $(-\Delta_t H-k^2 H=0)$  when  $x \in  \Omega_2\setminus \bar{\Omega}_1  $   $(t\in\Omega_2\setminus \bar{\Omega}_1)  $   respectively.
Since    $f$ vanishes on $\Omega_1\setminus \Omega$,
$$
H(x,t)=\int_{\Omega_1}  f(y) G_{\Omega_2}(x,y) G_{\Omega_2}(t,y) dy=\int_{\Omega}  f(y) G_{\Omega_2}(x,y) G_{\Omega_2}(t,y) dy.
$$
Thus, $H(x,t)=0$ follows from the assumption, for $t,x \in\Omega_2\setminus \bar{\Omega}.$   Note that $H$ also satisfies   $-(\Delta_x+\Delta_t) H-2k^2 H=0$, when $x,t \in  \Omega_2\setminus \bar{\Omega}_1  .$ The Unique Continuation Principle yields
$
H(x,t)=0~~\text{when~~}  t,x\in \bar{\Omega}_2\setminus  \Omega_1;
$
i.e.,
$$
\int_{\Omega_1}  f(y) G_{\Omega_2}(x,y) G_{\Omega_2}(t,y) dy=0,~~\text{for~} t,x\in \bar{\Omega}_2\setminus  \Omega_1.
$$
Applying  Lemma \ref{loc3},
we obtain the   equality $\int_{\Omega_1}  f u v=0 ,$ where
 $u,v\in C^{\infty}(\bar{\Omega}_1)$   are solutions to $-\Delta u-k^2 u=0$ and vanish  on $\partial \Omega_1\cap \partial \Omega_2$.
Combining   Proposition \ref{loc1}, $f$ vanishes on a neighborhood of $\partial \Omega_1\setminus (\partial \Omega_1\cap \partial \Omega_2)$. This implies that $f$ vanishes on a bigger neighborhood $\Theta_{\e}(t^{\prime})$ $t^{\prime}>t$ of the curve. Hence $I$ is open. This completes the proof of Theorem \ref{ds1}.
\hfill $\Box$



\section{Applications}
In this section, we give   simultaneous recovery of the partial data inverse problems for a class of nonlinear Helmholtz Schr\"{o}dinger equations   as an application of Theorem \ref{ds1}.

 In some applications, the discontinuity of the medium extends to the cases of obstacles or cavities embedded in the medium. Mathematically, equations are satisfied in the medium $\Omega$ minus the obstacle or cavity region $D$ and subject to certain boundary conditions at $\partial D$. Determining an unknown obstacle or cavity goes back to Schiffer's work  and is a substantial topic in inverse scattering theory, for example for sound waves  modeled by the linear wave equation and electromagnetic waves modeled by the linear Maxwell's equations. Among many existing methods (such as the probe method and so on), one method, known as the enclosure method,  uses CGO-type solutions to determine obstacles,  cavities and inclusions. The enclosure method was first introduced by
      Ikehata in   \cite{Ike99,Ike00} to reconstruct   a cavity $D$  inside a conductive medium $\Om$ for Schr\"{o}dinger operator  when the surrounding potential is known a priori. See  \cite{Isa90,Isa09}   for more uniqueness results of obstacle problems when knowing the potential and \cite{LL17}  for
        a simultaneous reconstruction of both  embedded obstacle and its surrounding potential.

\subsection{Simultaneous recovery of cavity and coefficients}

\begin{theo}\label{wx1}
   Let $\Omega$ be a bounded open set in $\R^n$,   $n\geq 2$ with smooth boundary $\partial \Omega$. Let $D_1,D_2\subset\subset \Omega$ be nonempty subsets with smooth boundaries such that $\Omega\setminus \bar{D}_j$ are connected.
   For $j=1,2$ let
   $q_j(x)\in C^{\infty}(\Omega\setminus \bar{D}_j)$  and consider the following nonlinear Helmholtz equation
\begin{align}\label{app1}
\left\{\begin{array}{ll}
  -\Delta u_j-k^2 u_j+q_j(x)u_j^2 =0  ~~& \text{in~}\Omega\setminus \bar{D}_j,\\
   u_j =0 ~~ &\text{on~} \partial D_j,\\
   u_j =f ~~&\text{on~} \partial \Omega.
   \end{array}
   \right.
\end{align}
Assume that the Dirichlet-Neumann map $\Lambda_{q_j}^{D_j}$
 satisfy that
\begin{align*}
  \Lambda_{q_1}^{D_1}(f) =\Lambda_{q_2}^{D_2}(f) \text{on~}\partial \Omega,~~\forall f\in C^s(\partial \Omega) \text{ ~is sufficiently small, }
\end{align*}
with   $s>1$ with $s\not\in\Z$. Then
 $D_1=D_2=D$ and
 $q_1=q_2$ in $\Omega\setminus\bar{D}.$
\end{theo}

\Proof We apply   the standard high order linearization method and and take the same notations  as shown in section \ref{yy1}.
Let  $f=\epsilon_1f_1+\epsilon_2f_2$ for the boundary value of \eqref{app1}, where $f_1$, $f_2\in C^s(\partial \Omega),$ and  $\epsilon_1,\epsilon_2$ are sufficiently small numbers.
%
Thus,  $v_j^{(l)}:= \partial_{ \epsilon_l}u_j\big|_{\epsilon_l=0}$ satisfy
\begin{align}\label{ap2vjl1}
\left\{\begin{array}{ll}
  -\Delta v_j^{(l)}-k^2 v_j^{(l)}    =0& \text{in~}\Omega \setminus \bar{D}_j,\\
   v_j^{(l)} = f_l&    \text{on~} \partial D_j,\\
  v_j^{(l)} =0& \text{on~} \partial \Omega,
   \end{array}
   \right.
\end{align}
with $j=1,2$ and $l=1,2.$

 Let $G$ be the connected component of $\Omega\setminus \bar{D}_1\cup \bar{D}_2$ of which the boundary contains $\partial \Omega.$
Let  $\tilde{v}^{(l)}=v_1^{(l)}-v_2^{(l)}$ for each $l$ in $G$. Then  $\tilde{v}^{(l)}$ satisfies
\begin{align*}
\left\{\begin{array}{ll}
  -\Delta \tilde{v}^{(l)}-k^2 \tilde{v}^{(l)}   =0&~\text{in~}G,\\
   \tilde{v}^{(l)}=\partial_{\nu}    \tilde{v}^{(l)}=0  & \text{on~} \partial \Omega,
   \end{array}
   \right.
\end{align*}
with $l=1,2.$ The unique continuation principle   of solutions to Helmholtz equation yields
$\tilde{v}^{(l)}=0$ in $G$. Therefore,
for each $l,$ $v_1^{(l)}=v_2^{(l)}$ in $G$. Hence $v_2^{(l)}=0$ on $\partial D_1\setminus D_2$ from the continuity.

To recovery the unknown boundary we need to  choose $f_l $.
Let $U=D_1\setminus \bar{D}_2$. Then
$\partial U=(\partial D_1\setminus D_2)\cup (\partial D_2\cap D_1)$.
Let $\varphi\in H^1(U)$ such that
 \begin{align*}
\left\{\begin{array}{ll}
  -\Delta \varphi-k^2 \varphi   =0&~\text{in~}U,\\
  \phi =0  & \text{on~} \partial D_2\cap D_1,\\
   \phi =g  & \text{on~}  \partial D_1\setminus D_2,\\
   \end{array}
   \right.
\end{align*}
where $g\in H^1(U)$   which has nonvanishing  trace   supported in  $\partial D_1\setminus D_2$.

Apply  Proposition (\ref{pen2}) with $U=D_1\setminus \bar{D}_2$ and $M=\Omega\setminus D_2$. Then $M\setminus \bar{U}=\Omega\setminus \bar{D}_1\cup \bar{D}_2$ is connected. Let  $\Gamma= \partial \Omega\subset \partial M$.  Hence, there exists some  $f_l$ with $\supp f_l\subset \Gamma$ and a corresponding solution $v_2^{(l)}$ to
\begin{align*}
\left\{\begin{array}{ll}
  -\Delta v_2^{(l)}-k^2 v_2^{(l)}   =0&~\text{in~}M,\\
  v_2^{(l)} =f_l  & \text{on~} \partial M,
   \end{array}
   \right.
\end{align*}
such that
  $v_2^{(l)}$ is close to $\phi$ in $H^1(U)$ and so is the trace of $v_2^{(l)}$ and $\phi$ on $\partial U$.
  Note that $v_2^{(l)} $ should vanish  on $\partial D_1\setminus D_2$. This contradicts with the assumption that the trace of  $\phi$ has nonvanishing support $g $ on $\partial D_1\setminus D_2$. This shows that $D_1=D_2.$
 From now on,
denote by $D=D_1=D_2$. It follows that
\begin{align}\label{ap2vjl3}
 v^{(l)}:= v_1^{(l)}=v_2^{(l)} ~\text{in~}\Omega\setminus \bar{D},
\end{align}
for $ l=1,2;$
and $ v^{(l)}$ is a solution to (\ref{ap2vjl1}).

To recover the coefficient,
 we choose $v^{(0)}$  to be a solution to
\begin{align}\label{ap2v0}
\left\{\begin{array}{ll}
  -\Delta v^{(0)}-k^2 v^{(0)}   =0&~\text{in~}\Omega\setminus \bar{D},\\
    v^{(0)}= 0  & \text{on~} \partial D,\\
  v_j^{(0)}=f_0&\text{on~} \partial \Omega,
   \end{array}
   \right.
\end{align}
where $f_0$ is not identically zero on $\partial \Omega.$

  Denoted by $w_j = \partial^2_{\epsilon_1   \epsilon_2}u_j\big|_{\e_1=\e_2=0}.$ It follows from the second order linearization that  $w_j = 0$ on $\pd D\cup \pd \Om$.
Moreover,  the DN map yields
$$
\partial_{\nu}w_1 |_{\partial \Omega}=\partial_{\nu}w_2 |_{\partial \Omega}.
$$
Thus, similar as in section \ref{yy1},
 one has the integral identity
$$
\int_{\Omega\setminus \bar{D}}  (q_1(x)-q_2(x))v^{(1)}v^{(2)} v^{(0)}dx=0,
$$
where  $v^{(1)}, v^{(2)}$  satisfy (\ref{ap2vjl1}) and  $ v^{(0)}$ satisfies (\ref{ap2v0}).
Apply Theorem \ref{ds1} and one has
\beq\label{iw2}
 (q_1(x)-q_2(x))v^{(0)}(x)=0
\eeq
identically in $\Omega\setminus \bar{D}$.
Similarly,
for any $\tilde{x}_0\in \Omega\setminus \bar{D}$, we   apply Proposition (\ref{pen2}) and construct $v^{(0)}$ a solution to (\ref{ap2v0}) by choosing some $f_0$, with $v^{(0)}(\tilde{x}_0)\not=0.$ Hence $(q_1-q_2)(\tilde{x}_0)=0$. Since $\tilde{x}_0$ is arbitrary in $\Omega \setminus \bar{D}$, then $q_1=q_2 $ in $\Omega\setminus \bar{D}$.
\hfill $\Box$

\subsection{Simultaneous recovery of boundary and coefficients}

\begin{theo}
   Let $\Omega_j$ be a bounded open set in $\R^n$,   $n\geq 2$ with smooth boundary $\partial \Omega$ for $j=1,2$. Consider $\Gamma\subset \partial \Omega$ be a proper nonempty open subset of $\partial \Omega_1\cap \partial \Omega_2.$
   For $j=1,2$ let
   $q_j(x)\in C^{\infty}(\Omega_j)$  and consider the following nonlinear Helmholtz equation
\begin{align}\label{app2}
\left\{\begin{array}{ll}
  -\Delta u_j-k^2 u_j+q_j(x)u_j^2=0& \text{in~}\Omega_j,\\
   u_j=0& \text{on~} \partial \Omega_j\setminus \Gamma,\\
   u_j=f&\text{on~}  \Gamma.
   \end{array}
   \right.
\end{align}
Assume that the Dirichlet-Neumann map $\Lambda_{q_j}^{\Omega_j}$
 satisfy that
\begin{align*}
  \Lambda_{q_1}^{\Omega_1}(f) =\Lambda_{q_2}^{\Omega_2}(f) \text{on~}\Gamma,~~\forall f\in C^s_c(\Gamma) \text{ ~is sufficiently small, }
\end{align*}
with   $s>1$ with $s\not\in\Z$. Then
 $\Omega_1=\Omega_2=\Omega$ and
 $q_1=q_2$ in $\Omega .$
\end{theo}

\Proof
 We apply   the standard high order linearization method and the procedures and  notations as in the proof of Theorem \ref{wx1}.
Let $f=\epsilon_1f_1+\epsilon_2f_2$.
Let  $v_j^{(l)}=\frac{\partial }{\partial \epsilon_l}\big|_{\e_l=0}u_j$ are solutions to the linearized Helmholtz equation with boundary value $f_l\in C_c^s(\Gamma),$ for $l=1,2$ and $j=1,2.$
  Similarly as in the proof of  Theorem \ref{wx1},  the unique continuation principle gives rise to $ v_1^{(l)}=v_2^{(l)}$ in $G$ for each $l $, where   $G$ is a connected component of $\Omega_1\cap \Omega_2$ of which the boundary contains $\Gamma.$  Hence $v_2^{(l)}=0$ on $\partial \Omega_1\cap \Omega_2$ from the continuity.


To recovery the unknown boundary we need to  choose the boundary values  $f_l $ for $l=1,2$.
Let $U=\Omega_2\setminus \bar{\Omega}_1$. Then
$\partial U=(\partial \Omega_2\setminus \Omega_1)\cup (\partial \Omega_1\cap \Omega_2)$.
Let $\varphi\in H^1(U)$ such that
 \begin{align*}
\left\{\begin{array}{ll}
  -\Delta \varphi-k^2 \varphi   =0&~\text{in~}U,\\
  \phi =0  & \text{on~} \partial \Omega_2\setminus \Omega_1,\\
   \phi =g  & \text{on~} \partial \Omega_1\cap\Omega_2,\\
   \end{array}
   \right.
\end{align*}
where $g\in H^1(U)$   which has nonvanishing  trace   supported in  $\partial \Omega_1\cap\Omega_2$.

Apply Proposition (\ref{pen2}) with $U=\Omega_2\setminus \bar{\Omega}_1$ and $M=\Omega_2$. Then $M\setminus \bar{U}=\Omega_1\cap \Omega_2$ is connected with $\Gamma=\Gamma_0\subset \partial M$.  Hence, there exists some  $f_l$ with $\supp f_l\subset \Gamma$ and a corresponding solution $v_2^{(l)}$ to
\begin{align*}
\left\{\begin{array}{ll}
  -\Delta v_2^{(l)}-k^2 v_2^{(l)}   =0&~\text{in~}M,\\
  v_2^{(l)} =f_l  & \text{on~} \partial M
   \end{array}
   \right.
\end{align*}
such that
  $v_2^{(l)}$ is close to $\phi$ in $H^1(U)$ and so is the trace of $v_2^{(l)}$ and $\phi$ on $\partial U$.
  Note that $v_2^{(l)} $ should vanish  on $\partial \Omega_1\cap \Omega_2$. This contradicts with the assumption that the trace of  $\phi$ has nonvanishing support $g $ on $ \partial \Omega_1\cap\Omega_2$. This shows that $\Omega_1=\Omega_2.$

 From now on,
denote by $\Omega=\Omega_1=\Omega_2$. It follows that
$
 v^{(l)}:= v_1^{(l)}=v_2^{(l)} ~\text{in~}\Omega,
$
for $ l=1,2;$
and $ v^{(l)}$ are solutions to the linearized Helmholtz equation with boundary value $f_l$ supported on $\Gamma $ for $ l=1,2.$
The recovery of the coefficient is then  similar to the proof of Theorem \ref{wx1}. By  taking the  second order linearization and   choosing some $v^{(0)}$ as in (\ref{ap2v0}) with $\pd D$ replaced by $\pd \Gamma$, one is able  to  construct an integral identity similar to (\ref{iw2}).  An application of  Theorem \ref{ds1} and similar arguments complete  the proof.
\hfill $\Box$
\newpage
\section{Nonlinear Maxwell equations}
In this section,  $\Om\su \R^3$ is a bounded domain with smooth boundary. The proof of Theorem \ref{KL1} and Theorem \ref{HL1} consist  of three parts. We first adopt  the standard procedures of high order  linearization of the admittance map to derive an integral identity.   Similarly,  we only provide key procedures of higher order linearization here and  refer readers to \cite{LLLS1,LLLS2} for more details. Then we follow the literature \cite{DKSU} to give  the uniqueness result locally in a small neighborhood of some point on the boundary. Lastly, we extend the local result to the global situation.
\subsection{Proof of Theorem \label{KL1}}

%

\Proof
 Let  $f=\e_1 f_{1}+\e_2 f_{2}+\e_3 f_{3}.$  $f_{l}\in C^{\infty}(\pd \Om)$ with $\supp f_{l} \su \Gamma_1$, for $l=1,2,3 $. Let $\e=(\e_1,\e_2,\e_3).$
 Denote
 $$W_j^{(l)}= \pd_{\e_l} H_j  |_{\e=0},~V_j^{(l)}= \pd_{\e_l} E_j  |_{\e=0},~W_j^{(123)}=\pd_{\e}H_j|_{\e=0},~ V_j^{(123)}=\pd_{\e}E_j|_{\e=0},$$
  for $j=1,2$ and $l=1,2,3 .$
Similarly as in section \ref{yy1}, one can show  that $ V^{(l)}:= V_1^{(l)}=   V_2^{(l)} $ and  $W^{(l)}:= W_1^{(l)} =   W_2^{(l)}    $ satisfy
\begin{align}
     \left\{
  \begin{aligned}
  \dt  V^{(l)}-i \omega \mu_0  W^{(l)}&=0,  \\
 \dt   W^{(l)}+i \omega \varepsilon_0  V^{(l)} &= 0,
  \end{aligned}
           \right.\label{firlin1}
\end{align}
in $\Om$ with $
 \nt V^{(l)} |_{\pd \Om}= f_{l} $ on the boundary $\pd \Om$, for $l=1,2,3 .$
By direct computation one can show that $(E,H)$    of     the following form is a  solution to \eqref{firlin1}
 \begin{align}
 \label{xz1}
 \begin{split}
   E&=\be_0^{-\frac12}\exi\left(
   \frac{1}{h^2} (\xi\cdot a)\xi-\frac{k}{h}\xi\times b-k^2 a
   \right), \\
   H&=\mu_0^{-\frac12}\exi\left(~\frac{1}{h^2} (\xi\cdot b)\xi+\frac{k}{h}\xi\times a-k^2 b \right),
   \end{split}
 \end{align}
  where $a,b\in \C^3 $ are complex vectors to be chosen later; and
  $\xi \in \C^3$ satisfies $\xi\cdot \xi=k^2h^2.$

The highest order linearization of the admittance map  yields
\begin{align*}
     \left\{
  \begin{aligned}
  \dt  V_j^{(123)}-i \omega \mu_0  W_j^{(123)}&= i \omega p_j  \sum_{(i_1i_2i_3)=\sigma(123)}W_j^{(i_1)} \cdot \overline{W}_j^{(i_2)} ~W_j^{(i_3)},\\
\dt  W_j^{(123)}-i \omega \be_0  V_j^{(123)}&= i \omega q_j  \sum_{(i_1i_2i_3)=\sigma(123)}V_j^{(i_1)} \cdot \overline{V}_j^{(i_2)}~V_j^{(i_3)},
  \end{aligned}
           \right.
\end{align*}
 in $\Om$ with  $\nt V_j^{(123)}|_{\pd \Om} = 0 $. Here, $\sigma\in S_3$ denotes  the permutation of the indices.
Thus,  we obtain the integral identity
 \begin{align}
   \int_{\Omega}& (p_1-p_2)\sum_{(i_1i_2i_3)=\sigma(123)}W^{(i_1)} \cdot \overline{W}^{(i_2)} ~W^{(i_3)}\cdot    \overline{W}^{(0)} dV,\label{inc1}
   \\
    &+\int_{\Omega} (q_1-q_2)\sum_{(i_1i_2i_3)=\sigma(123)}V^{(i_1)} \cdot \overline{V}^{(i_2)} ~V^{(i_3)}\cdot    \overline{V}^{(0)} dV=0,\label{inc2}
 \end{align}
by choosing $   W^{(0)}$ and $   V^{(0)}$ satisfying the following conjugate Maxwell's equations
\begin{align}
     \left\{
  \begin{aligned}
  \dt  V^{(0)}-i \omega \mu_0  W^{(0)}&=0, \\
 \dt   W^{(0)} +i \omega \varepsilon_0  V^{(0)} &= 0,
  \end{aligned}
           \right.\label{adj3}
\end{align}
in $\Om$ with boundary value $ \nt V^{(0)}|_{\pd \Om} = f_{0}$.
Here $f_{(0)}\in C^{\infty}(\pd \Om)$ with $\supp f_{0} \su \Gamma_2$.

%
%

\subsubsection{Local result}\label{secML}
Below we are to prove a local version of Theorem \ref{KL1}, i.e., to show $p_1=p_2$ and $q_1=q_2$ locally in a small neighborhood of some $x_0\in \Gamma_1$.  Without loss of generality, as in section 3, the domain
 can be set below the coordinate $\{x_1\leq 0\}$ by orthogonal transformation and translation, and $x_0=0\in \Ga_1$.


\begin{prop}\label{locmk1}
  Let $\Omega$ be a bounded open set in $\R^3$  with smooth boundary $\partial \Omega$. Let $x_0\in \Gamma_1\cap \Gamma_2\subset \partial \Omega$ be a convex point on the boundary such that $\Gamma\subset \Gamma_1\cap \Gamma_2$.
  Let $\Gamma$ be  a small neighborhood of $x_0$ on the boundary. Suppose that the cancellation (\ref{inc1}) and (\ref{inc2}) hold  for any
 smooth  solutions $(V^{(l)},W^{(l)})$  with $l=1,2,3$, and $(V^{(0)},W^{( 0)})$ to the systems  (\ref{firlin1})  and (\ref{adj3}), respectively. Then there exists $\delta>0$ such that $p_1-p_2 $ and $q_1-q_2 $ vanish  on $B(x_0,\delta)\cap \Omega$.
\end{prop}

\Proof
Let $p^{-1}(kh)=\{\xi\in C^3|\xi\cdot \xi=k^2h^2\}$, $\xi_0=(i,\sqrt{k^2h^2+1},0)\in p^{-1}(kh) $ and
$\eta_0=(i,-\sqrt{k^2h^2+1},0)\in p^{-1}(kh) $. Suppose $|\xi-\xi_0|<\e, \xi\in p^{-1}(kh)$, and $|\eta-\eta_0|<\e, \eta\in p^{-1}(kh)$. Then $\im(\xi_1)>0$ and $\im(\eta_1)>0$.

Below we construct solutions to  the systems  (\ref{firlin1})  and (\ref{adj3}), achieved by adding correction terms to the solutions satisfying the corresponding homogeneous Maxwell's equations   of the form \eqref{xz1}. Choose   $b=0$ in (\ref{xz1}). Denote
\beq
\begin{split}\label{choo1}
  E_1=\be_0^{-\frac12}\left( \frac{1}{h^2}(\xi\cdot a_1)\xi-k^2a_1 \right) \exi,~~~H_1=\mu_0^{-\frac12}  \frac{k}{h}(\xi\times  a_1)  \exi,\\
    E_2=\be_0^{-\frac12}\left( \frac{1}{h^2}(\xi\cdot a_2)\xi-k^2a_2 \right) \exi,~~~H_2=\mu_0^{-\frac12}  \frac{k}{h}(\xi\times  a_2)  \exi, \\
      E_3=\be_0^{-\frac12}\left( \frac{1}{h^2}(\xi\cdot a_3)\xi-k^2a_3 \right) \exi,~~~H_3=\mu_0^{-\frac12}  \frac{k}{h}(\xi\times  a_3)  \exi,
\end{split}
\eeq
and
\begin{align}
  E_0=\be_0 ^{-\frac12}\left( \frac{1}{h^2} (\bar{\eta}\cdot \bar{a}_0) \bar{\eta}-k^2 \bar{a}_0 \right) e^{\frac{i}{h}x\cdot \bar{\eta}},~~~H_0=-\mu_0^{-\frac12}  \frac{k}{h}(\bar{\eta} \times \bar{a}_0  ) e^{\frac{i}{h}x\cdot \bar{\eta}},\label{choo2}
\end{align}
where the constant vectors are defined by
\begin{align}\label{vir}
  a_1=(1,0,1),~~a_2=a_3=\bar{\xi}\times \bar{a}_1,~~a_0=\xi\times a_1.
\end{align}


Take a cutoff function  $0\leq \chi\leq 1\in C_c^{\infty}$   which equals 1 on $x_1\leq -2c$ containing $\partial \Omega\setminus \Gamma_1$ and $\partial \Omega\setminus \Gamma_2$
  and is compactly supported in $x_1\leq -c$. There exists
$(r_1, t_1)$ to the Dirichlet problem
\begin{align*}
     \left\{
  \begin{aligned}
  \dt  r_1-i \omega \mu_0  t_1&=0,\\
 \dt   t_1 +i \omega \varepsilon_0 r_1 &= 0,
  \end{aligned}
           \right.
\end{align*}
in $\Om$ with boundary value $ \nt r_1 |_{\pd \Om} = - \chi(x) E_1,$
satisfying the following estimates
\begin{align}\label{cr1}
\| r_1\|_{H^1 (\Omega)},  \| t_1\|_{H^1 (\Omega)}  \leq
 C(1+h^{-1}|\xi|)^{\frac32}e^{-\frac ch \im \xi_1 }e^{ \frac 1h |\im \xi^{'}| } ~~\text{if~} \im  \xi_1\geq 0.
\end{align}
Likewise, we have $(r_2, t_2)$, $(r_3, t_3)$  and $(r_0, t_0)$  with   similar estimates to (\ref{cr1}).
  Thus,
\begin{align*}
  V^{(1)} = E_1 + r_1, ~~~W^{(1)} = H_1 + t_1;~~~ V^{(2)} = E_2 + r_2,~~ W^{(2)} = H_2 + t_2;\\
    V^{(3)} = E_3 + r_3, ~~~W^{(3)} = H_3 + t_3;~~~ V^{(0)} = E_0 + r_0,~~ W^{(0)} = H_0 + t_0,
\end{align*}
are solutions to the systems (\ref{firlin1})  and (\ref{adj3}), respectively.
 Now, the derived  cancellation (\ref{inc1}) and (\ref{inc2})  become
 \begin{align}
   \int_{\Omega}& (p_1-p_2)\sum_{(i_1i_2i_3)=\sigma(123)}(H_{ i_1 }+t_{i_1}) \cdot (\overline{H_{ i_2 }+t_{i_2}})  (H_{ i_3 }+t_{i_3}) \cdot (\overline{H_{ 0 }+t_{0}})  dV\label{intenew1}
   \\
    &+\int_{\Omega} (q_1-q_2)\sum_{(i_1i_2i_3)=\sigma(123)}(E_{ i_1 }+r_{i_1}) \cdot (\overline{E_{ i_2 }+r_{i_2}})  (E_{ i_3 }+r_{i_3}) \cdot (\overline{E_{ 0 }+r_{0}})  dV=0.\label{intenew2}
 \end{align}
By virtue of (\ref{vir}), one has
$
H_1\cdot \bar{H}_2= H_1\cdot \bar{H}_3=H_1\cdot \bar{H}_0=0.
$
Therefore, only the items involving some $t_i$ survive in the expansion of the product
(\ref{intenew1}), which are exponentially decaying in view of (\ref{cr1}) when $\im \xi_1\geq 0$ and  $\im \xi^{\prime}=0$.

Next we calculate the items in the expansion of the product (\ref{intenew2}). Note that the
items involving some $r_i$ also have exponentially decaying properties from (\ref{cr1}) when
$\im \xi_1\geq 0$ and  $\im \xi^{'}=0$. We only need to calculate
$
(E_1\cdot \bar{E}_2)(E_3\cdot \bar{E}_0),~~(E_1\cdot \bar{E}_3)(E_2\cdot \bar{E}_0),~~(E_2\cdot \bar{E}_3)(E_1\cdot \bar{E}_0). $
Let $\xi=\xi_0+l_1$ and $\eta=\eta_0+l_2$, where $|l_1|,|l_2|\leq \e.$
\begin{small}
\begin{align*}
  & h^8 e^{ \frac{i}{h}x\cdot (2\xi-\bar{\xi} +\eta)} (E_1\cdot \bar{E}_2)(E_3\cdot \bar{E}_0)
  =    \bigg(
(\xi_0\cdot a_1)\xi_0-h^2k^2 a_1+ ((l_1\cdot a_1)\xi+ (\xi_0\cdot a_1) l_1)
  \bigg) \\
  &~~~~   \cdot
  \bigg(
   (\bar{\xi}_0\cdot \bar{a}_2)\bar{\xi}_0    -h^2k^2 \bar{a}_2
    +   \overline{((l_1\cdot a_2)\xi+(\xi_0\cdot a_2) l_1)}
  \bigg)
  \bigg(
  (\xi_0\cdot a_3)\xi_0-h^2 k^2 a_3
   + \\
   &~~ ~~~~~~((l_1\cdot a_3)\xi+
    (\xi_0\cdot a_3) l_1)
  \bigg)
     \cdot
     \bigg(
  (\eta_0\cdot a_0)\eta_0-h^2k^2 a_0+
 ( (l_2\cdot a_2 ) \eta+(\eta_0\cdot a_0 )l_2)
  \bigg)\\
\geq &  \big(
   (\xi_0\cdot a_1) (\bar{\xi}_0\cdot \bar{a}_2) (\xi_0\cdot \bar{\xi}_0)  +O(\e)  +O(h^2)  \big)
    \cdot
 \big(
  (\xi_0\cdot a_3) (\eta_0\cdot a_0) (\xi_0\cdot \eta_0)  +O(\e)+O(h^2)
  \big)\\
  \geq &  \big(
  4 +24\e |a_1| |a_2| |\xi_0|^3  +O(h^2) \big) \big(
  8 +2^5\e  |a_2|^2 |\xi_0|^3  +O(h^2) \big).
\end{align*}
\end{small}
Therefore, $h^8 e^{ \frac{i}{h}x\cdot (2\xi-\bar{\xi} +\eta)}(E_1\cdot \bar{E}_2)(E_3\cdot \bar{E}_0)$ is greater than some positive constant, provided $\e$ and $h$ are sufficiently small.
Direct calculations yield similar estimates for
  $(E_1\cdot \bar{E}_3)(E_2\cdot \bar{E}_0)$,  $(E_3\cdot \bar{E}_1)(E_2\cdot \bar{E}_0)$, $ (E_2\cdot \bar{E}_3)(E_1\cdot \bar{E}_0)$ and $(E_2\cdot \bar{E}_3)(E_1\cdot \bar{E}_0)$.
Hence,  when $\e$ and $h$ are sufficiently small, one has
$$
\left|\int_{  \Omega  } (q_1-q_2) \sum_{(i_1i_2i_3)=\sigma(123)}E_{i_1}\cdot \overline{E_{i_2}}~  E_{i_3}\cdot \overline{E_{i_0}} dV \right|\geq \left|\int_{  \Omega  } (q_1-q_2)h^{-8}  e^{-\frac ih x \cdot (2\xi-\bar{\xi}+\eta) } dV\right| \cdot 4$$
and the R.H.S. of the inequality is controlled by items involving at least one $t_{i}$ or $r_i$.
Combining the estimate (\ref{cr1}) yields
\begin{align*}
  \left|\int_{  \Omega  } (q_1-q_2)  e^{-\frac ih x \cdot (2\xi-\bar{\xi}+\eta) } dV\right|
\leq     C h^8\|  q_1-q_2  \|_{L^{\infty}(\Omega)}& (1+h^{-1}|\xi|)^{\frac92}(1+h^{-1}|\eta|)^{\frac32}\\
& e^{-\frac ch \min\{\im \xi_1, \im \eta_1\} }
e^{ \frac 1h( 3| \im \xi^{\prime}|+ |\im \eta^{\prime}|) },
\end{align*}
for $\im \xi_1\geq 0$ and $\im \eta_1\geq 0$. In particular, let
$$
z=2\xi-\bar{\xi}+\eta,~~z_0=(4i,0,0).
$$
When
$|\xi-a\zeta_0|<Ca\varepsilon $ and $|\eta-a\eta_0|<Ca\varepsilon$ with $\varepsilon\leq \frac{1}{2C}$, then $|z-az_0|<4Ca\varepsilon $, and
\begin{align*}
  \left|\int_{  \Omega  }  (q_1-q_2 ) e^{-\frac ih x \cdot z } dx\right|
\leq     C h^{2} \| q_1-q_2 \|_{L^{\infty}(\Omega)}
e^{-\frac {ca}{2h}   }
e^{ \frac {4Ca \varepsilon}{h} }.
\end{align*}
Conversely, for any  $z\in \C^3$, $|z-az_0|<4Ca\varepsilon $, $z$   can be decomposed to
$$
z=2\xi-\bar{\xi}+\eta,~~\xi, \eta\in p^{-1}(kh),~~|\xi-a\zeta_0|<Ca\varepsilon, |\eta-a\eta_0|<Ca\varepsilon.
$$
Therefore,   for all $z\in \C^n$, $|z-az_0|<4Ca \varepsilon$, the estimate holds
\begin{align*}
   \left|\int_{  \Omega  } (q_1-q_2)   e^{-\frac ih x \cdot z }  dx\right|
\leq   C h^{-3} \| q\|_{L^{\infty}(\Omega)}
e^{-\frac {ca}{2h}   }
e^{ \frac {2Ca \varepsilon}{h} }.
\end{align*}
Hence we are able to obtain
 an exponential decay for  the F.B.I. transform of $T(q_1-q_2) $  in  the  neighborhood  $|z-az_0|<4Ca \varepsilon$. Then  apply the water melon approach similarly to show  the uniqueness of $q(x)$ in a small neighborhood of $x_0=0$, i.e., there exists some $\delta_1>0$ such that
 $$
q_1(x)=q_2(x),~~\forall x\in \Omega,~~-\delta_1\leq x_1\leq 0.
 $$

Now let $a=0 $ in (\ref{xz1}) and choose $b_1,b_2,b_3,b_0$ similar as in  (\ref{vir}). It is easy to see that the expressions of $E$ and $H$ interchange in (\ref{choo1})-(\ref{choo2})  and the similar procedures follow. Hence in the same way we arrive at the uniqueness of $p(x)$ in a small neighborhood of $x_0=0$, i.e., there exists some $\delta_2>0$ such that
 $$
p_1(x)=p_2(x),~~\forall x\in \Omega,~~-\delta_2\leq x_1\leq 0.
 $$
This completes the proof of Proposition \ref{locmk1}.
\hfill $\Box$

\subsubsection{From local to global}\label{secMG}

To extend to global results we need some preparatory work.   Lemma \ref{rung} describes a density result involving solutions to Maxwell's equation as Lemma \ref{rung}.  Remark \ref{funda1} and  Remark \ref{funda2} provide some properties of Green functions of Maxwell's equations and Helmholtz equation, which will be needed in the sequel.

\begin{lem}\label{rung} Let $\Omega_1\subset\Omega_2\subset \R^3$ with smooth boundary and $\partial \Omega_1\cap \partial \Omega_2\not=\emptyset$.
Let $\mathcal{R}$ be the set of $  \{(v|_{\Omega_1}, w|_{\Omega_1})\in C^{\infty}(\Omega_1, \C^3)^2 \}$ such that
$$\left\{
\begin{array}{rl}
     \dt v-i \omega \mu_0 w&=a,\\
 \dt w+i \omega \varepsilon_0 v &= \bar{a},
    \end{array}
\right.~\text{in~} \Omega_2,
$$
with boundary value $ \nt v|_{\pd \Om_2} =0,$ for some
$a \in C_c^{\infty}(\Omega_2, \C^3)$ with $\supp a \subset \bar{\Omega}_2\setminus \Omega_1$. Let
$\mathcal{S}$ be the set of $  \{(\tilde{v} , \tilde{w} )\in \C^{\infty}(\Omega_1, \C^3)^2 \}$ such that
$$\left\{
\begin{array}{rl}
     \dt \tilde{v}-i \omega \mu_0 \tilde{w}&=0,\\
 \dt \tilde{w}+i \omega \varepsilon_0 \tilde{v} &= 0,
    \end{array}
\right.~\text{in~} \Omega_1,
$$
with boundary value  $ \nt \tilde{v}|_{\pd \Om_1\cap \pd \Om_2} =0.$
Then $\mathcal{R}$
is dense in the space
$\mathcal{S} $ with respect to the
$L^2(\Omega_1,\C^3)^2$ topology.
\end{lem}

\Proof
  Let $g=(g_1,g_2)\in L^2(\Omega_1,\C^3)^2$ and $g$ is orthogonal to $\mathcal{R}.$
It suffices to show that   $g$ is also orthogonal to $\mathcal{S}.$

Let $\tilde{g}_1$ and $\tilde{g}_2$ be the $L^2(\Omega_2, \C^3)$   extension of $g_1$ and $g_2$ by taking value 0 in $\Omega_2\setminus \Omega_1$.
 Suppose  $(f,h)$ are the solutions to the following nonhomogeneous Maxwell's equations
  $$\left\{
\begin{array}{rl}
     \dt f+i \omega \mu_0 h&=\tilde{g}_2,\\
 \dt h- i \omega \varepsilon_0 f &= \tilde{g}_1
    \end{array}
\right.~\text{in~} \Omega_2,
$$
with boundary value $  \nt f|_{\pd \Om_2} =0.$
Hence, for every $(v,w)\in \mathcal{R}$,
  \begin{align}
  &\lag  g_1, v\rag_{\Omega_1}+ \lag g_2, w\rag_{\Omega_1}
     = \lag \dt h- i \omega \varepsilon_0 f, v\rag_{\Omega_2}+\lag  \dt f+i \omega \mu_0 h, w\rag_{\Omega_2}\nonumber\\
     =&\int_{\Omega_2} h\cdot \overline{\dt v}dV+ \int_{\partial \Omega_2} \nt h\cdot \overline{ v} dS- \int_{\Omega_2} i \omega \varepsilon_0 f\cdot \overline{ v} dV \nonumber\\
    &~~~~+\int_{\Omega_2} f\cdot \overline{\dt w}dV+ \int_{\partial \Omega_2} \nt f\cdot \overline{ w} dS+ \int_{\Omega_2} i \omega \mu_0 h\cdot \overline{ w} dV\nonumber\\
     =&\int_{\Omega_2} h\cdot \overline{ \dt v-i \omega \mu_0 w  }dV+ \int_{\Omega_2} f\cdot \overline{ \dt w+ i \omega \varepsilon_0 v  }dV\nonumber\\
      =& \int_{\Omega_2} h\cdot \overline{ a  } dV+ \int_{\Omega_2} f\cdot a dV\label{balance}
      \end{align}
vanishes.  Suppose  $a = c(x)+d(x)i,$ with  $\supp c(x), d(x)\subset \bar{\Omega}_2\setminus \Omega_1$. Then (\ref{balance})   becomes
  \begin{align*}
    0&=\int_{\Omega_2} h\cdot  (c(x)-d(x)i) dV+ \int_{\Omega_2} f\cdot (c(x)+d(x)i)  dV\\
    &=\int_{\Omega_2}(  f+h)\cdot   c(x) dV+ \int_{\Omega_2} (f-h)\cdot  d(x)i   dV.
  \end{align*}
  Since $c(x)$ and $d(x)$ are arbitrary in $L^2(\Omega_2, \R^3)$,  $f+h $ and $f-h$ vanish on $\bar{\Omega}_2\setminus \Omega_1 $, which implies $f $ and $h$ vanish on $\bar{\Omega}_2\setminus \Omega_1.$
 Then
 \begin{align}\label{aubd1}
   \nt f|_{\partial \Omega_1\setminus \partial \Omega_2}=0~~\text{and}~~ \nt h|_{\partial \Omega_1\setminus \partial \Omega_2}=0.
 \end{align}
 Combining the assumption $\nt f|_{\partial \Omega_2  }=0$, one has
 \begin{align}\label{aubd2}
    \nt f|_{\partial \Omega_1}=0.
 \end{align}
  For any $(\tilde{v},\tilde{w})\in \mathcal{S}$,
  \begin{align*}
 &\lag  g_1, \tilde{v}\rag_{\Omega_1}+ \lag g_2, \tilde{w}\rag_{\Omega_1}
     = \lag \dt h- i \omega \varepsilon_0 f, \tilde{v}\rag_{\Omega_1}+\lag  \dt f+i \omega \mu_0 h, \tilde{w}\rag_{\Omega_1}\\
     =&\int_{\Omega_1} h\cdot \overline{ \dt \tilde{v}-i \omega \mu_0 \tilde{w}  }dV+ \int_{\Omega_1} f\cdot \overline{ \dt\tilde{ w}+ i \omega \varepsilon_0 \tilde{v}  }dV\\
    &~~~~+\int_{\partial \Omega_1} \nt h\cdot \overline{ \tilde{v}} dS+\int_{\partial \Omega_1} \nt f\cdot \overline{ \tilde{w}} dS,
      \end{align*}
  vanishes  by virtue of (\ref{aubd1}), (\ref{aubd2}) and $\nt \tilde{v}|_{ \partial \Omega_1\cap \partial \Omega_2}=0.$
   Therefore, $g$ is orthogonal to $\mathcal{S}.$
   \hfill $\Box$

\begin{rem}\label{funda1}
(i)The fundamental solution   $G^{\mu,\ep}= \left(
     \begin{array}{c}
       G^{E}  \\
        G^{H}
     \end{array}\right)=
      \left(\begin{array}{cc}
                   E^{(1)}   &  E^{(2)}   \\
                     H^{(1)}   &  H^{(2)}
                  \end{array}
                  \right)
$  consists of two parts $(E^{(1)}, H^{(1)} )$   and $(E^{(2)}, H^{(2)} )$  satisfying    the following Maxwell's equations, respectively:
 $$
  \left\{
  \begin{array}{rl}
     \dt  E^{(1)}(x,y)-i \omega \mu_0  H^{(1)}(x,y)&= \delta (x-y) I_3,\\
      \dt  H^{(1)}(x,y)+i \omega \varepsilon_0  E^{(1)}(x,y) &=0,
      \end{array}
      \right.
          $$
          in $\Om_2$ with boundary value $ \nt  E^{(1)}|_{\pd \Om_2} =0;$ and
          $$
        \left\{
  \begin{array}{rl}
          \dt  E^{(2)}(x,y)-i \omega \mu_0  H^{(2)}(x,y)&= 0,\\
 \dt  H^{(2)}(x,y)+i \omega \varepsilon_0  E^{(2)}(x,y) &= \delta (x-y) I_3,
       \end{array}
      \right.
$$
 in $\Om_2$ with boundary value $\nt  E^{(2)}|_{\pd \Om_2} =0$.

 (ii)Similarly, the fundamental solution
   $\tilde{G}^{\ep,\mu}= \left(
     \begin{array}{c}
       \tilde{G}^{E}  \\
        \tilde{G}^{H}
     \end{array}\right)=
   \left(\begin{array}{cc}
                   \tilde{E}^{(1)}   &  \tilde{E}^{(2)}   \\
                     \tilde{H}^{(1)}   &  \tilde{H}^{(2)}
                  \end{array}
                  \right)$
consists of  two parts which satisfy the Maxwell's equations
$$
  \left\{
  \begin{array}{rl}
     \dt  \tilde{E}^{(1)}(x,y)-i \omega \ep_0  \tilde{H}^{(1)}(x,y)&= \delta (x-y) I_3,\\
      \dt \tilde{ H}^{(1)}(x,y)+i \omega \mu_0  \tilde{E}^{(1)}(x,y) &=0,
      \end{array}
      \right.
      $$
      in $\Om_2$ with boundary value $  \nt  \tilde{H}^{(1)}|_{\pd \Om_2} =0 ;$ and
      $$
        \left\{
  \begin{array}{rl}
          \dt  \tilde{E}^{(2)}(x,y)-i \omega \ep_0  \tilde{H}^{(2)}(x,y)&= 0,\\
 \dt  \tilde{H}^{(2)}(x,y)+i \omega \mu_0  \tilde{E}^{(2)}(x,y) &= \delta (x-y) I_3,
      \end{array}
      \right.
$$
 in $\Om_2$ with boundary value $\nt  \tilde{H}^{(2)}|_{\pd \Om_2} =0.$
\end{rem}

\begin{rem}   \label{funda2}
 If $(v,w)$ satisfies the boundary value problem
   \begin{align*}
    \left\{
     \begin{array}{rl}
     \dt  v -i \omega \mu_0  w &= a,\\
 \dt  w +i \omega \varepsilon_0  v  &= \bar{a },
    \end{array}
  \right.~~~ \text{in~} \Omega_2,
  \end{align*}
  with boundary value  $\nt  v |_{\pd \Om_2} =0,$
 by direct computation    $(  \bar{w},\bar{v})$ is a solution to the conjugate Maxwell's equations
     \begin{align*}
     \left\{
     \begin{array}{rl}
         \dt  \bar{w} -i \omega \ep_0  \bar{v} &= a,\\
           \dt   \bar{v}+i \omega \mu_0  \bar{w}  &= \overline{a },
  \end{array}
  \right.~~~ \text{in~} \Omega_2,
  \end{align*}
 with boundary value  $\nt   \bar{v}|_{\pd \Om_2}  =0.$   This implies that the conjugate $\bar{w}$ plays a role of the electricity $E$ while the conjugate $\bar{v}$ works as the magnetism $H$.

\end{rem}

 Let $G_{\Omega_2}(x,y)$   be the matrix of  Green function of the   Helmholtz equation
 $
  (-\Delta_y-k^2)G_{\Omega_2}(x,y)=\delta(x-y) I_3, ~\text{in ~}  \Omega_2,
  $
  where $I_3$ is the $3\times 3 $ identity matrix.

\vspace{0.2 cm}
\noindent
\textit{Proof of Theorem \ref{KL1}}~ Below we show the extension of local result (Proposition \ref{locmk1}) to the global situation   following the same  process as in section 2.2.2.

Let $\bar{p}(y)=p_1-p_2$ and $\bar{q}(y)=q_1-q_2$.
Fix a point $x_1\in \Omega$ and let $\theta:[0,1]\ra \bar{\Omega}$ be a $C^1$ curve joining $x_0\in \Gamma$ to $x_1$ such that $\theta(0)=x_0$  and $\theta^{\prime}(0)$ is the interior normal to $\partial \Omega$ at $x_0$ and $\theta(t)\in \Omega$ for all $t\in (0,1].$  Suppose $\Theta_{\e}(t)$ is a small neighborhood of the curve $\theta([0,t])$ in $\Om$ and $I$ is the set of the   time where  $\bar{p}$ and $ \bar{q}$ vanish on  $\Theta_{\e}(t)$.
 Similarly, $I$ is a closed set and $I$ is also non-empty due to local results,  provided that $\e$ is small enough.
 It suffices to show that $I$ is open.

 Following the same setting as in section 2.2.2,  we have  $\partial \Theta_{\e} (t)\cap \partial \Omega\subset \Gamma,
 $ and
 $$
 \Omega\setminus \Theta_{\e} (t)\subset \Omega_1,~~\partial \Omega\setminus \Gamma\subset \partial \Omega\cap \partial \Omega_1;~~\partial \Omega\setminus \Gamma\subset \partial \Omega\cap \partial \Omega_1\subset  \partial \Omega\cap \partial \Omega_2.
 $$
Adopting the notation in Remark \ref{funda1}, let
$$   G_1(x_1,y)=  \left(
     \begin{array}{c}
       G_1^{E} (x_1,y) \\
        G_1^{H}(x_1,y)
     \end{array}\right),~
     \tilde{G}_1(x_1,y)  =
 \left(
     \begin{array}{c}
       \tilde{G}_1^{E} (x_1,y) \\
        \tilde{G}_1^{H}(x_1,y)
     \end{array}\right),$$  for $x_1,    y\in \Omega_2$.
Similarly, denote $  G_2(x_2,y), $ $ G_3(x_3,y) $ and $\tilde{G}_0  (x_0,y), $ $\tilde{G}_2  (x_2,y),$ $\tilde{G}_3  (x_3,y)$ for $x_2,x_3,x_0, y\in \Omega_2$, respectively.
%
When $x_1,x_2,x_3,x_0\in \Omega_2\setminus \Omega_1$,  define
\begin{align}
H(x_1,x_2,x_3,x_0)=&\int_{\Omega_1}  \bar{q}(y) K^E(x_1,x_2,x_3,x_0,y) +\bar{p}(y)K^H(x_1,x_2,x_3,x_0,y)dy\label{sign}\\
=&\int_{\Omega_1}    \bar{q}(y) \cdot G_{ 1}^E (x_1,y) \otimes   \tilde{G}_{  2}^H(x_2,y)  \cdot  G_{ 3}^E(x_3,y) \otimes   \tilde{G}_{ 0}^H (x_0,y) \nonumber
\\
&+  \bar{q}(y) \cdot G_{ 2}^E (x_2,y) \otimes   \tilde{G}^H_{  1}(x_1,y)  \cdot  G_{ 3}^E (x_3,y) \otimes   \tilde{G}^H_{ 0}(x_0,y)  \nonumber
\\
&+   \bar{q}(y)\cdot G_{ 1}^E (x_1,y) \otimes   \tilde{G}^H_{  3}(x_3,y)  \cdot  G_{ 2}^E (x_2,y) \otimes   \tilde{G}^H_{ 0} (x_0,y) \nonumber\\
&+ \bar{q}(y) \cdot G_{ 3}^E (x_3,y) \otimes   \tilde{G}_{  1}^H(x_1,y)  \cdot  G_{ 2}^E (x_2,y) \otimes   \tilde{G}^H_{ 0}(x_0,y) \nonumber\\
&+ \bar{q}(y) \cdot G_{ 2}^E (x_2,y) \otimes   \tilde{G}^H_{  3}(x_3,y)  \cdot  G^E_{ 1}(x_1,y) \otimes   \tilde{G}^H_{ 0}(x_0,y) \nonumber \\
&+\bar{q}(y) \cdot G_{ 3}^E(x_3,y)  \otimes   \tilde{G}^H_{  2}(x_2,y) \cdot  G_{ 1}^E (x_1,y) \otimes   \tilde{G}^H_{ 0} (x_0,y)\nonumber\\
& +
 \bar{p}(y) \cdot G_{ 1}^H (x_1,y) \otimes   \tilde{G}_{  2}^E(x_2,y)  \cdot  G_{ 3}^H(x_3,y) \otimes   \tilde{G}_{ 0}^E (x_0,y)\nonumber
\\
&+  \bar{p}(y) \cdot G_{ 2}^H (x_2,y) \otimes   \tilde{G}^E_{  1}(x_1,y)  \cdot  G_{ 3}^H(x_3,y) \otimes   \tilde{G}^E_{ 0}(x_0,y)\nonumber
\\
&+  \bar{p}(y)\cdot G_{ 1}^H (x_1,y) \otimes   \tilde{G}^E_{  3}(x_3,y)  \cdot  G_{ 2}^H (x_2,y) \otimes   \tilde{G}^E_{ 0} (x_0,y)\nonumber \\
&+ \bar{p}(y) \cdot G_{ 3}^H (x_3,y) \otimes   \tilde{G}_{  1}^E(x_1,y)  \cdot  G_{ 2}^H (x_2,y) \otimes   \tilde{G}^E_{ 0}(x_0,y) \nonumber\\
&+ \bar{p}(y) \cdot G_{ 2}^H (x_2,y) \otimes   \tilde{G}^E_{  3}(x_3,y)  \cdot  G^H_{ 1}(x_1,y) \otimes   \tilde{G}^E_{ 0}(x_0,y) \nonumber \\
&+ \bar{p}(y) \cdot G_{ 3}^H(x_3,y)  \otimes   \tilde{G}^E_{  2}(x_2,y) \cdot  G_{ 1}^H (x_1,y) \otimes   \tilde{G}^E_{ 0} (x_0,y)  \label{mafan1}
 dy.
\end{align}
Here  the  tensor sign $\otimes$ represents the tensor product of the matrices and is  also  consistent with the tensor of  distributions when   the matrices are viewed as distributions in the sense that
$$
\int(u_1(x_1)\otimes u_2(x_2))(\phi_1(x_1)\otimes \phi_2(x_2))dx_1 dx_2=\int u_1 \phi_1 dx_1 \int u_2 \phi_2 dx_2,
$$
for
 $u_1,u_2\in D^{\prime}(X_j)$ and $\phi_j\in C_0^{\infty}(X_j)$ where $X_j$ are open sets of $\R^3.$

 The term involving $\bar{p}(y)$
 plays the role of a magnetic field, while  the term  involving $\bar{q}(y)$ serves as an electric field in the coupled Maxwell's equations with vanishing tangential component on the boundary.
Moreover, for  $x_1,x_2,x_3,x_0\in \Omega_2\setminus \Omega_1$,
\begin{align*}
H(x_1,x_2,x_3,x_0)&=\int_{\Omega_1}  \bar{q}(y) K^E(x_1,x_2,x_3,x_0,y) +\bar{p}(y)K^H(x_1,x_2,x_3,x_0,y)  dy\\
&=\int_{\Omega}   \bar{q}(y) K^E(x_1,x_2,x_3,x_0,y) +\bar{p}(y)K^H(x_1,x_2,x_3,x_0,y)dy.
\end{align*}
 Note that   each entry of $H$ is a sum of terms involving $\bar{p}(y)$ and $\bar{q}(y)$ and satisfies the Helmholtz equation $  (-\Delta_y-k^2) u=0,  $  for $x_1,x_2,x_3,x_0 \in  \Omega_2\setminus \bar{\Omega}_1  .$
Therefore, for $x_1,x_2,x_3,x_0 \in\Omega_2\setminus \bar{\Omega},$  from the assumption we have $$\int_{\Omega }   \bar{q}(y) K^E(x_1,x_2,x_3,x_0,y) +\bar{p}(y)K^H(x_1,x_2,x_3,x_0,y)dy=0.$$
It follows that  $H(x_1,x_2,x_3,x_0)=0$ when $x_1,x_2,x_3,x_0\in \Omega_2\setminus \bar{\Omega}$.

 When  $x_1,x_2,x_3,x_0 \in  \Omega_2\setminus \bar{\Omega}_1  ,$ the entries of $H$  satisfy the Helmholtz equation $  (-\Delta_y-k^2) u=0.
 $
  Hence, it follows from  the Unique Continuation Principle that
$$
H(x_1,x_2,x_3,x_0)=0~~\text{when~~}  x_1,x_2,x_3,x_0\in \bar{\Omega}_2\setminus  \Omega_1;
$$
i.e., $~~\text{for~} x_1,x_2,x_3,x_0\in \bar{\Omega}_2\setminus  \Omega_1,$
\begin{align*}
\int_{\Omega_1}   \bar{q}(y) K^E(x_1,x_2,x_3,x_0,y) +\bar{p}(y)K^H(x_1,x_2,x_3,x_0,y) dy=0.
\end{align*}
Let   the left hand side of    (\ref{mafan1})   act on
$
(a_1 (x_1)\otimes  a_2 (x_2), a_3 (x_3)\otimes  a_4 (x_4)  )
$
for  $x_1,x_2,x_3,x_0 \in  \Omega_2 $. Here all
$ a_i$ are assumed to have  the forms
$
 a_i(x_i)=(l_i(x_i), \overline{l_i(x_i)})^T,
$
where each $ l_i $  is a smooth vector function  on $C^{\infty}(\Omega_2,\C^3)$ with support  in $\bar{\Omega}_2\setminus \Omega_1$.
For each $i $, let
\begin{align*}
 v_i(y) = \int_{\Omega_2} G^E_{i } (x_{i },y)
 a_{i }(x_i) dx_{i },~~w_i(y)=\int_{\Omega_2} G^H_{i } (x_{i },y)  a_{i }(x_i)  dx_{i },\\
\tilde{w}_i(y) = \int_{\Omega_2} \tilde{G}^E_{i } (x_{i },y)  a_{i }(x_i) dx_{i },~~\tilde{v}_i(y)=\int_{\Omega_2} \tilde{G}^H_{i } (x_{i },y)  a_{i }(x_i)  dx_{i }.
\end{align*}
 Then (\ref{mafan1}) becomes
\begin{align}
0=&\int_{\Omega_1}  \bar{q}(y) K^E(x_1,x_2,x_3,x_0,y) \vec{a} +\bar{p}(y)K^H(x_1,x_2,x_3,x_0,y)\vec{a} dy\nonumber\\
=&\int_{\Omega_1}
 \bar{q}(y)
(v_1  \cdot \tilde{v}_{  2} )   \cdot  (v_{ 3}  \cdot \tilde{v}_{ 0} )
+  \bar{q}(y)   (v_{ 2}    \cdot \tilde{v}_{  1}) \cdot  (v_{ 3}  \cdot  \tilde{v}_{ 0}) \nonumber
\\
&+   \bar{q}(y)
 (v_1  \cdot \tilde{v}_{  3} )   \cdot  (v_{ 2}  \cdot \tilde{v}_{ 0} )
+  \bar{q}(y)   (v_{ 3}    \cdot \tilde{v}_{  1}) \cdot  (v_{ 2}  \cdot  \tilde{v}_{ 0})
  \nonumber\\
&+ \bar{q}(y)
 (v_{ 2}    \cdot \tilde{v}_{  3}) \cdot  (v_{ 1}  \cdot  \tilde{v}_{ 0}) +  \bar{q}(y) (v_{ 3}    \cdot \tilde{v}_{  2}) \cdot  (v_{ 1}  \cdot  \tilde{v}_{ 0})
 \nonumber\\
& +
 \bar{p}(y)
(w_1  \cdot \tilde{w}_{  2} )   \cdot  (w_{ 3}  \cdot \tilde{w}_{ 0} )
+  \bar{p}(y)   (w_{ 2}    \cdot \tilde{w}_{  1}) \cdot  (w_{ 3}  \cdot  \tilde{w}_{ 0}) \nonumber
\\
&+   \bar{p}(y)
 (w_1  \cdot \tilde{w}_{  3} )   \cdot  (w_{ 2}  \cdot \tilde{w}_{ 0} )
+  \bar{p}(y)   (w_{ 3}    \cdot \tilde{w}_{  1}) \cdot  (w_{ 2}  \cdot  \tilde{w}_{ 0})
  \nonumber\\
&+ \bar{p}(y)
 (w_{ 2}    \cdot \tilde{w}_{  3}) \cdot  (w_{ 1}  \cdot  \tilde{w}_{ 0}) +  \bar{p}(y) (w_{ 3}    \cdot \tilde{w}_{  2}) \cdot  (w_{ 1}  \cdot  \tilde{w}_{ 0})
  \label{mafan2}
 dy
\end{align}

The fundamental solutions of Remark \ref{funda1} imply that
  $(v_i,w_i)$ is a  solution  to the inhomogeneous Maxwell's systems
\begin{align} \label{value1}
\left\{
\begin{array}{rl}
     \dt v_i-i \omega \mu_0 w_i&=l_i,\\
 \dt w_i+i \omega \varepsilon_0 v_i &=\overline{ l_i },
    \end{array}
\right.
\end{align}
in $\Om_2$ with boundary value $  \nt v_i|_{\pd \Om_2}=0,$
 while $(\tilde{v}_i,\tilde{w}_i)$ satisfies
$$\left\{
\begin{array}{rl}
     \dt \tilde{w}_i-i \omega \ep_0 \tilde{v}_i&= l_i,\\
 \dt \tilde{v}_i+i \omega \mu_0 \tilde{w}_i &=\overline{ l_i},
    \end{array}
\right.
$$
 in $\Om_2$ with boundary value $ \nt \tilde{v}_i|_{\pd \Om_2} =0 .$ Hence,  the uniqueness of solutions  gives
 $
 \tilde{v}_i=\bar{v}_i $ and $ \tilde{w}_i=\bar{w}_i $ in $ \Omega_2,$ for $   i=0,1,2,3.  $
Therefore, (\ref{mafan2}) yields
 \begin{align}
     0 &=    \int_{\Omega_1}  (p_1-p_2)\sum_{(i_1i_2i_3)=\sigma(123)}(w_{ i_1 } \cdot \overline{w}_{ i_2 })  (w_{ i_3 }\cdot    \overline{w}_{0}) dy\nonumber
   \\
    &\quad+\int_{\Omega_1} (q_1-q_2)\sum_{(i_1i_2i_3)=\sigma(123)}(v_{ i_1 } \cdot \overline{v}_{ i_2 })  (v_{ i_3 }\cdot    \overline{v}_{0}) dy,\label{sj1}
 \end{align}
 where $(v_i,w_i)$ satisfies   $(\ref{value1})$. Since
 all $ a_i  \in C^{\infty}(\Omega_2,\C^3)$ are   supported in $\bar{\Omega}_2\setminus \Omega_1$, the set of all the solutions  $(v_i,w_i)$ forms the set $\mathcal{R}$.
  Lemma \ref{rung}  states   that $\mathcal{R}$ is dense in the space $\mathcal{S} $ with respect to  the
$L^2(\Omega_1,\C^3)^2$ topology.
Thus for any $V_i$ and $W_i$ satisfying   inhomogeneous Maxwell's systems on $\Omega_1$ with $\nt V_i$ vanishing on  $\partial \Omega_1\cap \partial \Omega_2$, we can choose $v_i^{(n)}$ and $w_i^{(n)}$ converge to $V_i$ and $W_i$ in $L^2(\Omega_1,\C^3)$ respectively.

Note that
each $v_i^{(n)}$, $w_i^{(n)}$ and $V_i,W_i$ are smooth in $\Omega_2$. Then $\|v_i^{(n)}\|_{L^{\infty}(\Omega_2)}$ and $\|w_i^{(n)}\|_{L^{\infty}(\Omega_2)}$  are bounded for each fixed $i$.
Thus assume
$$
\|v_i^{(n)}\|_{L^{\infty}(\Omega_2)},\|w_i^{(n)}\|_{L^{\infty}(\Omega_2)},
\|V_i  \|_{L^{\infty}(\Omega_2)},\|W_i \|_{L^{\infty}(\Omega_2)}
\leq M.
$$
It is easy to  show that $(v_1^{(n)}\cdot \bar{v}_2^{(n)})( v_3^{(n)}\cdot \bar{v}_0^{(n)})  $ converges to $(V_1 \cdot \bar{V}_2)~(  V_3 \cdot \bar{V}_0)$  in $L^1(\Om_1)$   as $n\ra \infty$.
The  equality (\ref{sj1}) extends to
\begin{align*}
     0 &=    \int_{\Omega_1}  (p_1-p_2)\sum_{(i_1i_2i_3)=\sigma(123)}(W_{ i_1 } \cdot \overline{W}_{ i_2  } )(W_{ i_3 }\cdot    \overline{W}_{0}) dV\nonumber
   \\
    &\quad+\int_{\Omega_1} (q_1-q_2)\sum_{(i_1i_2i_3)=\sigma(123)}(V_{ i_1 } \cdot \overline{V}_{ i_2 })(V_{ i_3 }\cdot    \overline{V}_{ 0 }) dV,
 \end{align*}
where
 $(V_i,W_i)\in C^{\infty}(\bar{\Omega}_1,\C^3)^2$  is a solution  to   inhomogeneous Maxwell's systems and vanishes  on $\partial \Omega_1\cap \partial \Omega_2$.
Applying Proposition \ref{locmk1},   $\bar{p} $ and $\bar{q}$ vanish on a neighborhood of $\partial \Omega_1\setminus (\partial \Omega_1\cap \partial \Omega_2)$. Then $\bar{p} $ and $\bar{q}$ vanish  on a bigger neighborhood $\Theta_{\e}(t^{\prime})$ $t^{\prime}>t$ of the curve. Hence $I$ is open.
\hfill $\Box$
\vspace{5mm}

\subsection{Proof of Theorem \ref{HL1}}

 Below we  follow  the standard procedures of high order  linearization of the admittance map.
\Proof Let $f^{k\om} =\e_1 f^{k\om}_{1}+\e_2  f^{k\om}_{2}, $  $f_{l}^{k\om}\in C^{\infty}(\pd \Om)$ with $\supp f_{l} ^{k\om}\su \Gamma_1$, for $k=1,2$ and $l=1,2  $. Let $\e=(\e_1,\e_2)$.
Denote
$
 ( V^{k\om,j}_{l}, W^{k\om,j}_{l} )=(\pd \e_lE^{k\om,j}  \big|_{\e=0},   \pd \e_l H^{k\om,j} \big|_{\e=0})$.
  Similarly, one can show that
$
   V_{l}^{k\om }:=V_{l}^{k\om,1 }=   V_{l}^{ k\om,2}$ and $ W_{l}^{k\om}:= W_{l}^{k\om, 1}=   W_{l}^{k\om, 2}
$
 satisfy
\begin{align}\label{firse}
\left\{
\begin{array}{rl}
\dt V^{\om }_{l}-i \omega \mu_0 W^{\om }_{l}&=0,\\
 \dt W^{\om }_{l}+i \omega \varepsilon_0 V^{\om  }_{l} &= 0,\\
  \dt V^{2\om }_{l}-i2 \omega \mu_0 W^{2\om }_{l}&=0, \\
 \dt W^{2\om }_{l }+i 2\omega \varepsilon_0 V^{2\om }_{l} &=0,
 \end{array}
 \right.
\end{align}
 in $\Om$ with boundary values $ \nt V^{k\om }_{l}|_{\pd \Om } =f^{k\om}_{l} $ for $k=1,2$ and $l=1,2 .$ The second linearization yields that  $(W_{12 }^{k\om, j}, V_{12 }^{k\om, j} )  = (\partial^2_{\e}  H^{k\om, j}|_{\e=0}, \partial^2_{\e} E^{k\om, j} |_{\e=0})$
 satisfies
 \begin{align*}
\left\{
\begin{array}{rl}
\dt V_{12}^{\om,j}-i \omega \mu_0 W_{12}^{\om,j}&=0,\\
 \dt W_{12}^{\om,j}+i \omega \varepsilon_0 V_{12}^{\om,j} &=-i \omega \chi^{(2),j}\big( \overline{V_{1}^{\om }}\cdot V_{2}^{2\om }+
  \overline{V_{2}^{\om }}\cdot V_{1}^{2\om }
 \big),\\
  \dt V_{12}^{2\om,j}-i2 \omega \mu_0 W_{12}^{2\om,j}&=0, \\
 \dt W_{12}^{2\om,j}+i 2\omega \varepsilon_0 V_{12}^{2\om,j} &=-i   \omega \chi^{(2),j}    V^{\om }_{1 } \cdot V^{ \om }_{ 2},
 \end{array}
 \right.
\end{align*}
 in $\Om$ with boundary values
 $   \nt V_{12 }^{k\om,j}|_{\pd \Om} = 0$ with $k=1,2 $ and $    j=1,2$.
We are able to derive the following  integral identities
 \begin{align}
   \int_{\Omega}  (\chi^{(2),1}-\chi^{(2),2}) \cdot \overline{V_{0}^{\om}}~
   \big( \overline{V_{1}^{\om }}\cdot V_{2}^{2\om }+
  \overline{V_{2}^{\om }}\cdot V_{1}^{2\om }
 \big)
   dV=0,\label{inte1}
    \end{align}
    and
  \begin{align}
   \int_{\Omega}  (\chi^{(2),1}-\chi^{(2),2}) \cdot \overline{V_{0}^{2\om}}~\big(
 V_{1}^{\om } \cdot V_{2}^{ \om } \big)
    dV=0, \label{inte2}
 \end{align}
by choosing  $  (V_{0}^{k\om},W_{0}^{k\om})$ to be the solution to the following conjugate Maxwell's equations
\begin{align}
\left\{
\begin{array}{rl}
\dt V^{\om}_{0}-i \omega \mu_0 W^{\om}_{0}&=0,\\
 \dt W^{\om}_{0}+i \omega \varepsilon_0 V^{\om}_{0} &=0,\\
  \dt V^{2\om}_{0}-i2 \omega \mu_0 W^{2\om}_{0}&=0, \\
 \dt W^{2\om}_{0}+i 2\omega \varepsilon_0 V^{2\om}_{0} &=0,
 \end{array}
 \right.\label{adjze}
\end{align}
in $\Om$ with boundary values   $\nt V_{0}^{k\om}|_{\pd \Om } = f_{0}^{k\om},$ where $f^{k\om}_{0} \in C^{\infty}(\pd \Om)$ with supports on $\Gamma_2$, for each $k$.
Next we construct
$ V_{1}^{\om } $,  $ V_{2}^{ \om } $ and $V_{0}^{2\om}$    satisfying     (\ref{firse}) and (\ref{adjze}) to plug in the integral identity (\ref{inte2}).
 \begin{rem}\label{Maxsol} It is easy to check that
  $(E,H)$   of the following form
 \begin{align}
   E &=\be_0^{-\frac12}\exi \tilde{A}_{\xi},~~\text{satisfying~$   \tilde{A}_{\xi}\cdot \xi=0$,}\label{firdef}\\
     H &=\be_0^{-\frac12}\exi \tilde{B}_{\xi},~~\text{with~$   \tilde{B}_{\xi}=-\frac{1}{kh}\xi \times   \tilde{A}_{\xi} $,}\label{firdef22}
 \end{align}
   satisfies the linearized Maxwell equations  by noting that
  $$
 \tilde{A}_{\xi}= -\frac{1}{k^2h^2}\cdot(
  (\xi\cdot\xi)   \tilde{A}_{\xi}- (\xi \cdot   \tilde{A}_{\xi})\xi
  )= -\frac{1}{k^2h^2} \xi \times  (\xi \times   \tilde{A}_{\xi})= \frac{1}{kh}\xi \times   \tilde{B}_{\xi}.
  $$
\end{rem}

\subsubsection{Local result}\label{secHL}
Below we are to prove a local version of Theorem \ref{HL1}, i.e., to show   $\chi^{(2),1}=\chi^{(2),2}$ locally in a small neighborhood of some $x_0\in \Gamma_1$. Without loss of generality we assume that  the domain
is   below the coordinate $\{x_1\leq 0\}$ by orthogonal transformation and translation and $x_0=0$.

\begin{prop}\label{locmh1}
  Let $\Omega$ be a bounded open set in $\R^3$  with smooth boundary $\partial \Omega$. Let $x_0\in \Gamma_1\cap \Gamma_2\subset \partial \Omega$ be a convex point on the boundary such that $\Gamma\subset \Gamma_1\cap \Gamma_2$.
  Let $\Gamma$ be  a small neighborhood of $x_0$ on the boundary. Suppose that the cancellation (\ref{inte2}) hold  for any
 smooth  solutions $(V_{1}^{\omega},W_{1}^{\omega})$, $(V_{2}^{\omega},W_{2}^{\omega})$, and $(V_{0}^{2\omega},W_{0}^{2\omega})$  to the systems  (\ref{firse}) and (\ref{adjze}). Then there exists $\delta>0$ such that $\chi^{(2),1}=\chi^{(2),2}$   on $B(x_0,\delta)\cap \Omega$.
\end{prop}

\Proof
\noindent
 (i) Define the set $p^{-1}(kh)=\{\xi\in \C^3|\xi\cdot \xi=k^2h^2\}$ and choose $$  \xi_0  =(i,\sqrt{1+k^2h^2 },0),~~ \eta_0=(i,0, \sqrt{1+k^2h^2 } )\in p^{-1}(kh).
$$
Suppose $|\xi-\xi_0|<\e, \xi\in p^{-1}(kh)$, and $|\eta-\eta_0|<\e, \eta\in p^{-1}(kh)$. Then $\im(\xi_1)>0$ and $\im(\eta_1)>0$.
Below we construct pairs of solutions  as in subsection \ref{KL1}. Here we technically construct $(E_{0}^{2\omega},H_{0}^{2\omega})$ of the expressions (\ref{firdef}) (\ref{firdef22}) and $(E_{1}^{\omega},H_{1}^{\omega})$, $(E_{2}^{\omega},H_{2}^{\omega})$
  of the expression (\ref{xz1}).

As in  Remark \ref{Maxsol}, choose
$ E_0^{2\om} =\be_0^{-\frac12}e^{-\frac{i}{h} x\cdot \varsigma }\tilde{A}_{\varsigma},$ where
\begin{align}\label{chig1}
\begin{split}
\varsigma & =(\sqrt{2(1-k^2h^2)}i, \sqrt{ 1+k^2h^2 },\sqrt{ 1+k^2h^2 }  )\in p^{-1}(2kh),\\
   \tilde{A}_{\varsigma}& = \left(-\sqrt{2(1-k^2h^2)}i, -\frac{(1-k^2h^2)}{ \sqrt{ 1+k^2h^2 }},-\frac{(1-k^2h^2)}{ \sqrt{ 1+k^2h^2 }} \right),
\end{split}
\end{align}
Since $\varsigma \cdot \tilde{A}_{\varsigma}=0 $,  $E_0^{2\om}$ together with some determined $H_0^{2\om}$ are solutions to the linearized Maxwell's equations following Remark \ref{Maxsol}.

Let  $b=0$, $  a_1=a_2=(1,0,0) $  in (\ref{xz1}). Denote  $H_1^{\om}=\mu_0^{-\frac12}  \frac{k}{h}(\xi\times  a_1)  \exi$,  $H_2^{\om}=\mu_0^{-\frac12}  \frac{k}{h}(\eta\times  a_2)   e^{-\frac{i}{h} x\cdot \eta }$,  and
\beq
\begin{split}
  E_1^{\om}=\be_0^{-\frac12}\exi \tilde{A}_{\xi}=\be_0^{-\frac12}\left( \frac{1}{h^2}(\xi\cdot a_1)\xi-k^2a_1 \right) \exi, \label{chn1}\\
       E_2^{\om}=\be_0^{-\frac12} e^{-\frac{i}{h} x\cdot \eta }  \tilde{A}_{\eta}=\be_0^{-\frac12}\left( \frac{1}{h^2}(\eta\cdot a_2)\xi-k^2a_2 \right)  e^{-\frac{i}{h} x\cdot \eta },~~~
       \end{split}
       \eeq
When $|\xi-\xi_0|<\e,$
\begin{align*}
 \tilde{A}_{\xi}&= h^{-2}\left(  (\xi\cdot a_1)\xi-k^2h^2a_1\right) =  h^{-2}\xi\times (\xi\times a_1)\\
&= h^{-2}  \big(  \xi_0\times (\xi_0\times a_1)+O(\e)\big)\\
&=h^{-2}  \bigg(  (
- 1-k^2h^2, \sqrt{1+k^2h^2}i,0)
+O(\e)\bigg),
\end{align*}
where $O(\e)$ is some terms independent of $x$.
Similarly,
when  $|\eta-\eta_0|<\e,$
\begin{align*}
 \tilde{A}_{\eta}&= h^{-2}  \left( (\eta\cdot a_2)\eta-k^2h^2a_2 \right)
 =h^{-2}   \eta\times (\eta\times a_2)\\
&= h^{-2}    \big(  \eta_0\times (\eta_0\times a_2)+O(\e)\big)\\
&=h^{-2}  \bigg(  (
- 1-k^2h^2,0, \sqrt{1+k^2h^2}i )
+O(\e)\bigg),
\end{align*}
where $O(\e)$ is some terms independent of $x$.
Construct
$(r_1^{\om}, t_1^{\om})$, $(r_2^{\om}, t_2^{\om})$ and $(r_0^{2\om}, t_0^{2\om})$    with   exponentially decaying  estimates as in (\ref{cr1}) such that
\begin{align*}
  V_1^{\om} = E_1^{\om} + r_1^{\om}, ~~~W_1^{\om} = H_1^{\om} + t_1^{\om};&~~~ V_2^{\om} = E_2^{\om} + r_2^{\om},~~ W_2^{\om} = H_2^{\om} + t_2^{\om};\\
      V_0^{2\om}  = E_0^{2\om}  + r_0^{2\om},&~~ W_0^{2\om}  = H_0^{2\om}  + t_0^{2\om}.
\end{align*}
are solutions to (\ref{firse}) and (\ref{adjze}), respectively.
Thus, the integral identity  (\ref{inte2}) becomes
 \begin{align}
   \int_{\Omega}& (\chi^{(2),1}-\chi^{(2),2}) \cdot
   \overline{(E_0^{2\om}  + r_0^{2\om})} ~~
    (E_{ 1}^{\om}+r_{ 1}^{\om}) \cdot (E_{ 2}^{\om}+r_{ 2}^{\om}) dV=0.\label{inten2}
 \end{align}

Next we calculate the  terms in the expansion of the product (\ref{inten2}). Note that the
terms involving some $r_i$ have exponentially decaying property from (\ref{cr1}) when
$\im \xi_1\geq 0$ and  $\im \xi^{'}=0$. As a result,
 \begin{align}
   \int_{\Omega}  (\chi^{(2),1}-\chi^{(2),2}) \cdot
   \overline{ E_0^{2\om}}   ~~
    ( E_{ 1}^{\om} \cdot  E_{ 2}^{\om}  ) dV \text{~consists of terms involving some $r_i^{k\om}$}.\label{intenew3}
 \end{align}
 For the L.H.S. of (\ref{intenew3}), one has
\begin{align*}
&\left| \int_{\Omega}   \ep_0^{-\frac32} e^{-\frac{i}{h}x\cdot(\xi+\eta-\bar{\varsigma})} (\chi^{(2),1}-\chi^{(2),2}) \cdot
   \overline{ \tilde{A}_{\varsigma}}   ~~
    ( \tilde{A}_{\xi} \cdot  \tilde{A}_{\eta}  ) dV\right|\\
    =& \bigg| \bigg( (
- 1-k^2h^2, \sqrt{1+k^2h^2}i,0)\cdot (
- 1-k^2h^2,0,  \sqrt{1+k^2h^2}i )
+O(\e)\bigg)  \bigg| \\
  & \cdot
  \bigg|   \int_{\Omega}  h^{-4} \ep_0^{-\frac32} e^{-\frac{i}{h}x\cdot(\xi+\eta-\bar{\varsigma})} (\chi^{(2),1}-\chi^{(2),2}) \cdot
   \overline{ \tilde{A}_{\varsigma}}
     dV\bigg|\\
 \geq & C   h^{-4} \ep_0^{-\frac32}\bigg|   \int_{\Omega}   e^{-\frac{i}{h}x\cdot(\xi+\eta-\bar{\varsigma})} (\chi^{(2),1}-\chi^{(2),2}) \cdot
   \overline{ \tilde{A}_{\varsigma}}
     dV\bigg|,
\end{align*}
 when $\e$ and $h$ is sufficiently small.
Denote by
\beq\label{yyx}
f(x;h)=  (\chi^{(2),1}-\chi^{(2),2}) \cdot
   \overline{ \tilde{A}_{\varsigma}}.
\eeq
Hence, when $\im \xi_1\geq 0$ and $\im \eta_1\geq 0$,
\begin{align*}
\bigg|   \int_{\Omega}   e^{-\frac{i}{h}x\cdot(\xi+\eta-\bar{\varsigma})} f(x;h)
     dV\bigg|&\leq C h^4\ep_0^{\frac32}
       \| \chi^{(2),1}-\chi^{(2),2} \|_{L^{\infty}(\Omega)} (1+h^{-1}|\xi|)^{\frac32}(1+h^{-1}|\eta|)^{\frac32}\\
&(1+h^{-1}|\varsigma|)^{\frac32} e^{-\frac ch \min\{\im \xi_1, \im \eta_1\} }
e^{ \frac 1h( | \im \xi^{'}|+ |\im \eta^{'}|) }.
\end{align*}
 In particular, let
   $$z=\xi+\eta-\bar{\varsigma},~~z_0=((2+\sqrt{2-2k^2h^2})i,0,0).$$
   When
$|\xi-a\xi_0|<Ca\varepsilon $ and $|\eta-a\eta_0|<Ca\varepsilon$ with $\varepsilon\leq \frac{1}{2C}$,   $|z-az_0|< 2C a\varepsilon, $ and
\begin{align*}
  \left|\int_{  \Omega  }  f(x;h) e^{-\frac ih x \cdot z } dx\right|
\leq     C  h^{-\frac{1}{2}}\ep_0^{\frac32}
       \| \chi^{(2),1}-\chi^{(2),2} \|_{L^{\infty}(\Omega)}
e^{-\frac {ca}{2h}   }
e^{ \frac {4Ca \varepsilon}{h} }.
\end{align*}
Conversely, for any  $z\in \C^3$, $|z-az_0|<2Ca\varepsilon $, $z$   can be decomposed to
$$
z= \xi  +\eta- \bar{\varsigma},~~\xi, \eta\in p^{-1}(kh),~~|\xi-a\zeta_0|<Ca\varepsilon, |\eta-a\eta_0|<Ca\varepsilon.
$$
Thus, following the previous procedures,
one can show the exponentially decaying property of the F.B.I. transform $Tf$
for all $x\in \Omega,~|x_1|\leq \delta_1,$ with $\delta_1$ sufficiently small.
  Hence, in view of  (\ref{chig1}) (\ref{yyx}) and
 $$
 \lim_{h\ra 0}(2\pi h)^{-\frac{n}{2}}Tf(x)=(\chi^{(2),1}-\chi^{(2),2}) \cdot   \left( \sqrt{2 }i, - 1,-1 \right) ~~\text{in $L^p(\Om)$},
 $$
 one has
\begin{align}
  (\chi^{(2),1}-\chi^{(2),2}) \cdot   \left( \sqrt{2 }i, - 1,-1 \right)=0,~~x\in \Om, -\delta_1\leq x_1\leq 0.\label{conclu1}
\end{align}

\vspace{0.5 cm}
\noindent
 (ii)
Now  choose
$$\xi_0=(i,-\sqrt{1+k^2h^2 },0),~\eta_0=(i,0, -\sqrt{1+k^2h^2 } ) \in p^{-1}(kh). $$   Suppose $|\xi-\xi_0|<\e, \xi\in p^{-1}(kh)$, and $|\eta-\eta_0|<\e, \eta\in p^{-1}(kh)$. Then $\im(\xi_1)>0$ and $\im(\eta_1)>0$.

Let $ E_0^{2\om} =\be_0^{-\frac12}e^{-\frac{i}{h} x\cdot \varsigma }\tilde{A}_{\varsigma} $  as in  Remark \ref{Maxsol}, where
\begin{align}\label{bxd}
\begin{split}
\varsigma  &=(\sqrt{2(1-k^2h^2)}i,- \sqrt{ 1+k^2h^2 },-\sqrt{ 1+k^2h^2 }  )\in p^{-1}(2kh),\\
   \tilde{A}_{\varsigma}& = \left(-\sqrt{2(1-k^2h^2)}i,  \frac{(1-k^2h^2)}{ \sqrt{ 1+k^2h^2 }}, \frac{(1-k^2h^2)}{ \sqrt{ 1+k^2h^2 }}. \right)
\end{split}
\end{align}
Again $E_0^{2\om}$ with some determined $H_0^{2\om}$ are  solutions to the linearized Maxwell's equations    by noting
$\varsigma \cdot \tilde{A}_{\varsigma}=0 $. $ E_1^{ \om} $ and $ E_2^{ \om} $ are defined as in (i) of the form  (\ref{xz1}), where $b=0$, $  a_1=a_2=(1,0,0) $.
 Let
   $$z=\xi+\eta-\bar{\varsigma},
   ~~z_0=((2+\sqrt{2-2k^2h^2})i,0,0).$$
Similarly, the amplitudes of $ E_1^{ \om} $ and $ E_2^{ \om} $  have the estimates
\begin{align*}
 \tilde{A}_{\xi}& = h^{-2}  \big(  \xi_0\times (\xi_0\times a_1)+O(\e)\big)\\
&=h^{-2}  \bigg(  (
- 1-k^2h^2, -\sqrt{1+k^2h^2}i,0)
+O(\e)\bigg),
\end{align*}
and
\begin{align*}
 \tilde{A}_{\eta}&=   h^{-2}  \bigg(  (
- 1-k^2h^2,0, -\sqrt{1+k^2h^2}i )
+O(\e)\bigg),
\end{align*}
 where $O(\e)$ is some terms independent of $x$, when $|\xi-\xi_0|<\e$ and $|\eta-\eta_0|<\e,$ respectively.

After going through similar procedures as before, we  arrive at an exponential decaying property of $Tf$, where $f$ is defined in \eqref{yyx}.  Recall the expression \eqref{bxd} of $\tilde{A}_{\varsigma}$. Then,
  $$
 \lim_{h\ra 0}(2\pi h)^{-\frac{n}{2}}Tf(x)=(\chi^{(2),1}-\chi^{(2),2}) \cdot   \left( \sqrt{2 }i,   1, 1 \right) ~~\text{in $L^p(\Om)$}.
 $$
Therefore, one has
\begin{align}
  (\chi^{(2),1}-\chi^{(2),2}) \cdot   \left( \sqrt{2 }i,   1, 1 \right)=0,~~x\in \Om, -\delta_2\leq x_1\leq 0.\label{conclu2}
\end{align}

\vspace{0.5 cm}
\noindent
 (iii)
We choose
$$\xi_0=(i, \sqrt{1+k^2h^2 },0),~\eta_0=(i,0, -\sqrt{1+k^2h^2 } )\in p^{-1}(kh) .$$    Suppose $|\xi-\xi_0|<\e, \xi\in p^{-1}(kh)$, and $|\eta-\eta_0|<\e, \eta\in p^{-1}(kh)$. Then $\im(\xi_1)>0$ and $\im(\eta_1)>0$.

Let  $E_0^{2\om} =\be_0^{-\frac12}e^{-\frac{i}{h} x\cdot \varsigma }\tilde{A}_{\varsigma},$ where
\begin{align}\label{bxd2}
\begin{split}
\varsigma & =(\sqrt{2(1-k^2h^2)}i,  \sqrt{ 1+k^2h^2 },-\sqrt{ 1+k^2h^2 }  )\in p^{-1}(2kh),\\
   \tilde{A}_{\varsigma}& = \left(-\sqrt{2(1-k^2h^2)}i, - \frac{(1-k^2h^2)}{ \sqrt{ 1+k^2h^2 }}, \frac{(1-k^2h^2)}{ \sqrt{ 1+k^2h^2 }} \right),
\end{split}
\end{align}

Again $E_0^{2\om}$ together
with some determined $H_0^{2\om}$ are   solutions to the linearized Maxwell equations    by noting
$\varsigma \cdot \tilde{A}_{\varsigma}=0 $ and  Remark \ref{Maxsol}.

Similarly, take
$  a_1=a_2=(1,0,0).
$
When $|\xi-\xi_0|<\e,$
\begin{align*}
 \tilde{A}_{\xi}& = h^{-2}  \big(  \xi_0\times (\xi_0\times a_1)+O(\e)\big)\\
&=h^{-2}  \bigg(  (
- 1-k^2h^2,  \sqrt{1+k^2h^2}i,0)
+O(\e)\bigg),
\end{align*}
where $O(\e)$ is some terms independent of $x$;
and
when  $|\eta-\eta_0|<\e,$
\begin{align*}
 \tilde{A}_{\eta}&=   h^{-2}  \bigg(  (
- 1-k^2h^2,0, -\sqrt{1+k^2h^2}i )
+O(\e)\bigg),
\end{align*}
where $O(\e)$ is some terms independent of $x$. Let
   $$z=\xi+\eta-\bar{\varsigma},~~
   z_0=((2+\sqrt{2-2k^2h^2})i,0,0).$$
Following similar procedures as before and noting    the expressions  \eqref{yyx} and \eqref{bxd2} for $f$ and   $\tilde{A}_{\varsigma}$, we arrive at
  $$
 \lim_{h\ra 0}(2\pi h)^{-\frac{n}{2}}Tf(x)=(\chi^{(2),1}-\chi^{(2),2}) \cdot   \left( \sqrt{2 }i,   -1, 1 \right) ~~\text{in $L^p(\Om)$}.
 $$
It follows that
\begin{align}
  (\chi^{(2),1}-\chi^{(2),2}) \cdot   \left( \sqrt{2 }i,  - 1, 1 \right)=0,~~x\in \Om, -\delta_3\leq x_1\leq 0.\label{conclu3}
\end{align}
In summary of  (\ref{conclu1}) (\ref{conclu2}) (\ref{conclu3}), one has
  \begin{align*}
    \chi^{(2),1}-\chi^{(2),2}=0,~~x\in \Om, -\delta\leq x_1\leq 0,
  \end{align*}
for some positive $\delta=\min\{\delta_1,\delta_2,\delta_3\}.$
This completes the proof of Proposition \ref{secHL}.

\subsubsection{From local to global}
\label{secHG}
\vspace{0.2 cm}
\noindent
\textit{Proof of Theorem \ref{HL1}}~ Below we show the extension of local result (Proposition \ref{secHL}) to the global situation   following the same  process as in section 2.2.2 and  \ref{secMG}.

 Let $\bar{\chi}(y)=\chi^{(2),1}-\chi^{(2),2}$.  Fix a point $x_1\in \Omega$ and let $\theta(t) $ be a parameterized   curve in $\Om$ joining $x_0\in \Gamma$ to $x_1$ such that $\theta(0)=x_0$  and $\theta^{\prime}(0)$ is the interior normal to $\partial \Omega$ at $x_0$.    Suppose $\Theta_{\e}(t)$ is a small neighborhood of the curve $\theta([0,t])$ in $\Om$ and $I$ is the set of the   time where  $\bar{\chi} $ vanishes on  $\Theta_{\e}(t)$.    The set $\Omega\cup B(x_0,\e^{\prime})$ is smoothed out into an open set $\Omega_2$ with smooth boundary such that
$$
\partial \Omega\setminus \Gamma\subset \partial \Omega\cap \partial \Omega_1\subset  \partial \Omega\cap \partial \Omega_2,
$$  and it suffices to show the set $I$ is open.

Following the notations in Section \ref{secMG},
 define
\begin{align*}
H(x_1,x_2,x_0)
= \int_{\Omega_1}    \bar{\chi}(y) \otimes \tilde{G}^H_{ 0} (x_0,y) \cdot G_1^E(x_1,y) \otimes   G_{  2}^E (x_2,y)
 dy,
\end{align*}
for $x_1,x_2,x_0\in \Omega_2\setminus \Omega_1$.
Moreover,
\begin{align}
H(x_1,x_2,x_0)&=\int_{\Omega_1}  \bar{\chi}(y) \otimes \tilde{G}^H_{ 0} (x_0,y) \cdot G_1^E(x_1,y) \otimes   G_{  2}^E (x_2,y).
\end{align}
When $x_1,x_2,x_0 \in\Omega_2\setminus \bar{\Omega},$ from the assumption we have   $H(x_1,x_2,x_3,x_0)=0$, as  all the entries of $H$  plays the role of an electric field which satisfies the homogeneous Maxwell's equations with vanishing tangential components on the boundary  for $x_1,x_2,x_0 \in  \Omega_2\setminus \bar{\Omega}_1  .$

Note that
 all the entries of $H$  satisfy the Helmholtz equation $  (-\Delta_y-k^2) u=0,  $  for $x_1,x_2,x_0 \in  \Omega_2\setminus \bar{\Omega}_1  .$  Hence, from the unique continuation principle,
$$
H(x_1,x_2,x_0)=0~~\text{when~~}  x_1,x_2,x_0\in \bar{\Omega}_2\setminus  \Omega_1;
$$
i.e., $~~\text{for~} x_1,x_2,x_3,x_0\in \bar{\Omega}_2\setminus  \Omega_1,$
\begin{align}\label{final1}
\int_{\Omega_1}   \bar{\chi}(y)\otimes
\tilde{G}^H_{ 0} (x_0,y) \cdot G_1^E(x_1,y) \otimes   G_{  2}^E (x_2,y)dy=0.
\end{align}
Let the operator at the left hand side of    (\ref{final1})   act on
 $$
  \left(1\otimes
            a_0 (x_0),
            a_1 (x_1) \otimes
            a_2 (x_2) \right)
 $$
for  $ x_0,x_1,x_2\in  \Omega_2 $. Here all $ a_i $ are of the forms
$
a_i(x_i)=(l_i(x_i), \overline{l_i(x_i)})^T,
$
and all $l_i$ are smooth vector functions on $C^{\infty}(\Omega_2,\C^3)$ with supports contained in $\bar{\Omega}_2\setminus \Omega_1$.
Denote by
$$
 v_i(y) = \int_{\Omega_2} G^E_{i } (x_{i },y)
 a_{i } (x_{i }) dx_{i },
~~\text{and}~~
\tilde{v}_{0}(y)=\int_{\Omega_2} \tilde{G}^H_{0 } (x_{0 },y) a_{0 }(x_0)  dx_{0 },
$$ for $i=1,2$.
Thus
  $v_i$ and $w_i$ are solutions to the inhomogeneous Maxwell's systems (\ref{value1})
and similarly, one has
 $
 \tilde{v}_0=\bar{v}_0$ and $  \tilde{w}_0=\bar{w}_0$ in $ \Omega_2,$
 where $(v_0,w_0)$ satisfies   similar Maxwell's equations   (\ref{value1}) with frequency  replaced by $2\om$.
Therefore,   (\ref{final1}) becomes
 \begin{align*}
     0 &=    \int_{\Omega_1}  \bar{\chi} \cdot \bar{v}_{0}~~    v_{1} \cdot v_{2}  dy,
 \end{align*}
 where $(v_i,w_i)$ satisfy the system  $(\ref{value1})$. Since
 all $ a_i  \in C^{\infty}(\Omega_2,\C^3)$ with supports contained in $\bar{\Omega}_2\setminus \Omega_1$, the solutions  $(v_i,w_i)$ form the set $\mathcal{R}$.

  Lemma \ref{rung}  states   that $\mathcal{R}$ is dense in the space $\mathcal{S} $ with respect to  the
$L^2(\Omega_1,\C^3)^2$ topology. In the same  way as in Section \ref{secMG}, one can show that
the  equality (\ref{sj1}) extends to
\begin{align*}
     0 =    \int_{\Omega_1} \bar{\chi}\cdot \bar{V}_{0}~~    V_{1} \cdot V_{2}
     dV
 \end{align*}
 holds for  every
 $V_i,W_i\in C^{\infty}(\bar{\Omega}_1,\C^3)$  which satisfy   inhomogeneous Maxwell systems and of which tangential components vanish  on $\partial \Omega_1\cap \partial \Omega_2$.
Hence by virtue of Proposition \ref{loc1} the cancellation states that $\bar{\chi}$ vanishes on a neighborhood of $\partial \Omega_1\setminus (\partial \Omega_1\cap \partial \Omega_2)$. This implies that $\bar{\chi}$ vanishes on a bigger neighborhood $\Theta_{\e}(t^{\prime})$ $t^{\prime}>t$ of the curve. Hence $I$ is open.

\vspace{5mm}

\appendix
\section{Appendix}\label{sec:appdx}

\subsection{Well-posedness of nonlinear Helmholtz Schr\"{o}dinger equation }\label{Appendix:well-posedness}

Let $\Om \su \R^n$ be a bounded domain with smooth boundary. Consider the boundary value problem for the nonlinear Helmholtz equation with small boundary value
\begin{align}\label{wel1}
\left\{\begin{array}{rl}
  -\Delta u-k^2 u+q(x,u)=0& \text{in~}\Omega,\\
   u=f& \text{on~}  \partial \Om.
   \end{array}
   \right.
\end{align}
Here $q(x,u)$ is smooth in $u$ and $\pd_z^k q(x,u_0)\in L^{\infty}(\Om\times \R) $   for $k\in \N$.
Moreover, $q(x,0)=\pd_u q(x,0)=0.$

\begin{theo}\label{aw1}
  Suppose that $p>\frac n2$. There exists a discrete subset $\Sigma\su \R$ such that
 for every $k\not\in \Sigma$,    the  boundary value problem   (\ref{wel1}) has a unique solution    $u(x)$
satisfying
$$
\|u\|_{W^{2,p}(\Omega)}\leq c \|f\|_{W^{2-\frac1p,p}(\pd \Om)},
$$
whenever $\|f\|_{W^{2-\frac1p,p}(\pd \Om)}<\delta$ is sufficiently small. Here, $c>0$ is a positive constant  depending only  on $p,k,\Om  $ and the coefficients of $q.$
\end{theo}

 \begin{Proof}
  Given $F\in L^p(\Om)$, from Fredholm theory, there exists a discrete subset $\Sigma\su \R$ such that
 for any $k\not\in \Sigma$, the Dirichlet boundary problem
 \begin{align*}
\left\{\begin{array}{rl}
  -\Delta u-k^2 u =F& \text{in~}\Omega,\\
   u=0& \text{on~}  \partial \Om,
   \end{array}
   \right.
\end{align*}
 has a unique solution $u\in W_0^{2,p}(\Om)$ satisfying  $\|u\|_{ W_0^{2,p}(\Om)}\leq c \|F\|_{L^p(\Om)},$ the constant $c$ depending only  on $p,k$ and $\Om  $.
In other words, for any $k$ outside $\Sigma$,
 the operator $\mathcal{L}=-\Delta-k^2:   W^{2,p}_0( \Om)\ra L^p(\Om)$ has a continuous inverse.

As a result, one can  show that,  for any $k$ outside $\Sigma$ and   $f\in W^{2-\frac1p,p}(\pd\Om)$,
 the linearized equation
  \begin{align*}
\left\{\begin{array}{rl}
  -\Delta u-k^2 u  =0 &  \text{in~}\Omega,\\
   u =f&   \text{on~}  \partial \Om,
   \end{array}
   \right.
\end{align*}
of (\ref{wel1})  has a solution $u_0\in  W^{2,p}(\Om)$ satisfying $\|u_0\|_{ W^{2,p}(\Om)}\leq c \|f\|_{W^{2-\frac1p,p}(\pd\Om)} \leq c\delta$.
Hence, $u=u_0+v$ is a solution to the  boundary value problem   (\ref{wel1}), if $v$ satisfies
  \begin{align*}
\left\{\begin{array}{cl}
  -\Delta v-k^2 v = -q(x,u_0+v)  & \text{in~}\Omega,\\
   v= 0  & \text{on~}  \partial \Om.
   \end{array}
   \right.
\end{align*}

Denote by $\mathcal{F}(v)=-q(x,u_0+v)$.  Since $q(x,0)=\pd_u q(x,0)=0,$ then
$q(x,u)=q_r(x,u)u^2$, where $q_r(x,u)=\int_0^1\pd_u^2 q(x,tu)$ $(1-t) dt $ is bounded. Furthermore, $q_r(x,u)$ is also Lipschitz in $u$ noting $\pd^k_u q(x,u)\in L^{\infty}(\Om \times \C)$ for $k\in \N$.
Combining  that   $W^{2,p}(\Om) \hookrightarrow C(\bar{\Om})$ by Sobolev
embedding theorem for $p>\frac n2$,  for $v_1 $ and $v_2$ in
 $X_{\delta}$ defined as
$$X_{\ep}=\{v\in W^{2,p}(\Om)|~ \|v\|_{W^{2,p}(\Om)} \leq c\ep\},$$
one has
$\|v_1\|_{L^{\infty}(\Om)}, \|v_2\|_{L^{\infty}(\Om)}\leq c\ep,$ and $\|u\|_{L^{\infty}(\Om)}\leq c\delta,$.
Thus,
\begin{align*}
  &\|\mathcal{L}^{-1}\circ\mathcal{F}(v_1)-\mathcal{L}^{-1}\circ\mathcal{F}(v_2)\|_{W^{2,p}(\Om)}
\\
  \leq&\|\mathcal{F}(v_1)-\mathcal{F}(v_2) \|_{L^p(\Om)}
  \\ =&\|(q_r(x,u_0+v_1)-q_r(x,u_0+v_2))(u_0+v_2)^2\|_{L^p(\Om)}\\
  &+\|q_r(x,u_0+v_1)( (u_0+v_2)^2-(u_0+v_1)^2 )\|_{L^p(\Om)}\\
  \leq& c \|v_1-v_2 \|_{L^p(\Om)}\|u_0+v_2\|^2_{L^{\infty}(\Om)}+  c\|2u_0+v_1+v_2\|_{L^{\infty}(\Om)}\|v_1-v_2\|_{L^p(\Om)}\\
  \leq & c(\delta^2+\ep^2+\ep\delta+\ep+\delta) \|v_1-v_2 \|_{L^p(\Om)}.
\end{align*}
The constant is independent of $\ep$ and $\delta$.
This implies that $\mathcal{L}^{-1}\circ\mathcal{F}$ is a contraction on $X_{\delta}$ provided that $\ep$ and $\delta$ are sufficiently small. It follows that there exists some $v\in X_{\delta}$ such that
$v=\mathcal{L}^{-1}\circ\mathcal{F}(v)$.
Furthermore,
\begin{align*}
\|v\|_{W^{2,p}(\Om)}&\leq c\|\mathcal{F}(v)\|_{L^p(\Om)}\leq c (\|u_0\|^2_{L^{\infty}(\Om)}+\|v\|^2_{L^{\infty}(\Om)})\\
&\leq c(\delta \|f\|_{W^{2-\frac 1p,p}(\pd\Om)} +\ep\|v\|_{W^{2,p}(\Om)}).
\end{align*}
Thus one has $\|v\|_{W^{2,p}(\Om)} \leq   \|f\|_{W^{2-\frac 1p,p}(\pd\Om)} $ provided $\delta$ and $\ep$ sufficiently small.
In summary, $u=u_0+v\in W^{2,p}(\Om)$ is a solution to the nonlinear Helmholtz equation with small boundary value satisfying the required  estimate.
\hfill $\Box$

\end{Proof}

\subsection{Runge approximation}

We need the following  Runge approximation property in the previous sections.

Let $M\subset \R^n,n\geq 2 $ be a bounded open set with $C^{\infty}$ boundary.
Consider the dual space of $H^{1}(M)$, given by
$$
\tilde{H}^{-1}(M)=\left\{f \in  H^{-1 }(\R^n)|~\supp(f)\subset \bar{M}\right\}.
$$
For $f\in \tilde{H}^{-1 }(M) $ and $g\in H^{1 }(M)$, let $\tilde{g}=\mathrm{Ext}(g)\in H^{1 }(\R^n)$ be an extension of $g$.
%
%

Define
the following duality pairing
$$
( f,g)_{\tilde{H}^{-1 }(M),H^{1 }(M) } =( f, \tilde{g})_{\tilde{H}^{-1 }(\R^n),H^{1 }(\R^n) }=\int_{\R^n} f \overline{\tilde{g}}  dx.
$$
%
\begin{prop}\label{pen2}
Let $U \subset   M $ be a domain with $C^{\infty}$ boundary such that $M\setminus \bar{U}$ is connected. Assume that $\partial U\cap \partial M=\bar{V}$, where $V\subset \partial M$ is open with $C^{\infty}$ boundary. Let $\Gamma$ be a nonempty open subset of $\partial M$, and $\bar{\Gamma}\cap \bar{V}=\emptyset.$
 Then the set
$$
\mathcal{R}=\{u|_{U}:~ -\Delta u-k^2 u=0~\text{in}~M, u=g~\text{on}~\partial M; g\in C_c^{\infty}(\Gamma)\}
$$
is dense in the space
$\mathcal{S} =\{u\in H^{1 }(U):~-\Delta u-k^2u=0~\text{in}~U,u=0~\text{on}~  \partial U\cap \partial M\}$ with respect to the
$H^{1 }(U)$ norm.
\end{prop}

\Proof
  Let $f\in  \tilde{H}^{-1 }(U)$ be orthogonal to $\mathcal{R}$, i.e., one has
$
( f, u|_{U})_{\tilde{H}^{-1 }(U),H^{1 }(U) } =0
$ for any $u|_{U}$ satisfying
$$ -\Delta u-k^2 u=0~\text{in}~M, u=g~\text{on}~\partial M; g\in C_c^{\infty}(\Gamma).$$
From Hahn-Banach theorem it suffices to show that
$
( f, v)_{\tilde{H}^{-1 }(U),H^{1 }(U) }  =0
$ for any $v\in \mathcal{S}. $

Construct $M_0$ with smooth boundary such that $M\subset M_0$ and $ \partial M\setminus \bar{\Gamma}\subset \partial M\cap \partial M_0\ $. Thus, $\partial U   \cap \partial M\subset \partial U\cap \partial M_0$.
 For $f\in \tilde{H}^{-1 }(U) $, there exists a sequence $\{f_j\}\subset C_0^{\infty}(U) $ such that
 $f_j\ra f$ in $\tilde{H}^{-1 }(U) .$

 Let $w\in H_0^{1 }(M_0)$ and $w_j\in C^{\infty}(\bar{M_0})\cap H_0^{1 }(M_0) $ be the unique solutions to the following Helmholtz equations
\begin{align*}
\left\{
\begin{array}{ll}
  -\Delta w-k^2 w =f~&\text{in}~M_0, \\
  w =0~&\text{on}~\partial M_0,
  \end{array}
  \right.
\end{align*}
 and
 \begin{align*}
 \left\{
\begin{array}{ll}
  -\Delta w_j-k^2 w_j=f_j~&\text{in}~M_0, \\
  w_j=0~&\text{on}~\partial M_0.
    \end{array}
  \right.
\end{align*}
Since  $f_j\ra f$ in $\tilde{H}^{-1 }(U) ,$ the regularity of solutions to the elliptic equation yields
$w_j\ra w$ in $H^{1 }(M) .$

 Let $\varphi\in H^{1}(M) $ and $\varphi_j\in C^{\infty}(\bar{M})\cap H^{1}(M) $ be the unique solutions to the following Helmholtz equation
   \begin{align*}
   \left\{
\begin{array}{ll}
  -\Delta \varphi -k^2 \varphi =0~&\text{in}~M, \\
     \varphi =w ~&\text{on}~\partial M.
\end{array}
     \right.
\end{align*}
and
  \begin{align*}
     \left\{
\begin{array}{ll}
  -\Delta \varphi_j-k^2 \varphi_j=0~&\text{in}~M, \\
     \varphi_j=w_j~&\text{on}~\partial M.
     \end{array}
     \right.
\end{align*}
Note that $w_j\ra w$ in $H^{1}(M) ,$ then
 $\varphi_j\ra \varphi$ in $H^{1}(M) .$

 Let $\omega_j|_{M}=w_j-\varphi_j$ and $\omega|_{M}=w-\varphi $. Then $\omega$ and $\omega_j$ satisfy
  \begin{align} \label{bir}
        \left\{
\begin{array}{ll}
  -\Delta \omega -k^2 \omega =f~&\text{in}~M, \\
    w =0~&\text{on}~\partial M.
    \end{array}
     \right.
\end{align}
and
 \begin{align*}
      \left\{
\begin{array}{ll}
  -\Delta \omega_j-k^2 \omega_j=f_j~&\text{in}~M, \\
     w_j=0~&\text{on}~\partial M,
  \end{array}
     \right.
\end{align*}
 respectively.
 Thus, for $u\in \mathcal{R},$ one has
 \begin{align*}
   0&=( f, u|_{U})_{\tilde{H}^{-1 }(U),H^{1 }(U) } =\lim_{j\ra \infty}( f_j, u|_{U})_{\tilde{H}^{-1 }(U),H^{1 }(U) } =\lim_{j\ra \infty}\int_{M}  -\Delta w_j u-k^2 w_j  u dV\\
   &=\lim_{j\ra \infty} \int_{M} ( -\Delta u -k^2u) w_j    dV  + \lim_{j\ra \infty}\int_{\partial M}       \partial_{\nu} \omega_j  u - \omega_j \partial_{\nu}   u dS\\
   &=\lim_{j\ra \infty}\int_{\Gamma}       \partial_{\nu} \omega_j  g dS= \int_{\Gamma}       \partial_{\nu} \omega  g dS
 \end{align*}
 Hence, $   \partial_{\nu} \omega|_{\Gamma}=0.$
In summary, $\omega$ satisfies
  \begin{align*}
  \left\{
  \begin{array}{ll}
  -\Delta \omega -k^2 \omega =0~&\text{in}~M \setminus \bar{U}, \\
  \omega =0~&\text{on}~\partial M\\
  \partial_{\nu} \omega=0~&  \text{on}~\Gamma
    \end{array}
     \right.
\end{align*}
 The unique continuation principle implies that
 $$
 \omega|_{M\setminus \bar{U}}=0.
 $$
 Recall that $\omega\in H^{1}_0(M)$,then $\omega|_{\partial U\cap \partial M}=0$. Hence, $\omega|_{\partial U}=0.$

 Thus $\omega $ may be identified with an element in $H^{1}_0(U)$.
There exists $\psi_j\in C_0^{\infty}(U)$ such that $ \psi_j\ra \omega$ in $H^{1}(U)$. Hence $-\Delta \psi_j -k^2\psi_j\ra -\Delta \omega -k^2  \omega$ in $H^{-1}(\R^n).$

 For any $v\in \mathcal{S},$ let $\tilde{v}$ be the extension of $v$ to $\R^n.$
 \begin{align*}
&(-\Delta \omega -k^2  \omega, \tilde{v} )_{H^{-1}(\R^n), H^1(\R^n) }\\
=&\lim_{j\ra \infty} (-\Delta \psi_j -k^2\psi_j,\tilde{v}  )_{H^{-1}(\R^n), H^1(\R^n) }
 = \lim_{j\ra \infty} \int_{M} ( -\Delta \psi_j -k^2\psi_j) v    dV \\
      =& \lim_{j\ra \infty} \int_{U} ( -\Delta v -k^2v) \psi_j    dV+\lim_{j\ra \infty}\int_{\partial U}       \partial_{\nu} \psi_j  v - \psi_j \partial_{\nu}   v dS
       = 0.
 \end{align*}

Consider $\hat{f}=f-(-\Delta \omega -k^2  \omega ).$ $\hat{f}\in H^{-1}(\R^n) $ and  $H^{-\frac12}(\partial U)$ with $\supp \hat{f}\su \partial U$.
      From [*],
$$
\hat{f}=\tilde{f}\otimes \delta_{\partial U},~~\tilde{f}\in H^{-\frac12}(\partial U).
$$
Following  (\ref{bir}) and $f\in \tilde{H}^{-1}(U)$ one has $\supp \hat{f}\su \bar{V}=\partial U\cap \partial M.$ Thus $\supp \tilde{f}\su \bar{V}.$
Therefore,
there exists $\tilde{f}_j\in C_0^{\infty}(V)$ such that $\tilde{f}_j\ra \tilde{f}$ in $H^{-\frac12}(\partial U).$ Then
\begin{align*}
  (\hat{f},v|_{U})_{\tilde{H}^{-1 }(U),H^{1 }(U) } &=(\hat{f}, \tilde{v} )_{H^{-1}(\R^n), H^1(\R^n) }=( \tilde{f}, v|_{\partial U})_{H^{-\frac{1}{2}}(\partial U), H^{\frac12}(\partial U) }\\
  =& \lim_{j\ra \infty}( \tilde{f}_j, v|_{\partial U})_{H^{-\frac{1}{2}}(\partial U), H^{\frac12}(\partial U) }= \lim_{j\ra \infty}\int_{\partial U} \tilde{f}_j v dS=0,
\end{align*}
since $\supp \tilde{f}_j\su V$ and  $v|_{\bar{V}}=0.$ Therefore,
\begin{align*}
  (f,v|_{U})_{\tilde{H}^{-1}(U), H^1(U)  }&=(-\Delta \omega -k^2  \omega, \tilde{v} )_{H^{-1}(\R^n), H^1(\R^n) }+( \hat{f}, \tilde{v} )_{H^{-1}(\R^n), H^1(\R^n) }=0.
\end{align*}
This completes the proof of Proposition \ref{pen2}.
\hfill $\Box$

\end{document}